# The Square Sieve and a Lang-Trotter Question for Generic Abelian Varieties

Samuel Bloom*

February 16, 2017

**Abstract**

Let $A$ be a $g$-dimensional abelian variety over $\mathbb{Q}$ whose adelic Galois representation has open image in $\operatorname{GSp}_{2g}\widehat{\mathbb{Z}}$. We investigate the "Frobenius fields" $\mathbb{Q}(\pi_p) = \operatorname{End}(A_p) \otimes \mathbb{Q}$ of the reduction of $A$ modulo primes $p$ at which this reduction is ordinary and simple. We obtain conditional and unconditional asymptotic upper bounds on the number of primes at which $\mathbb{Q}(\pi_p)$ is a specified number field and, when $A$ is two-dimensional, at which $\mathbb{Q}(\pi_p)$ contains a specified real quadratic number field. These investigations continue the investigations of variants of the Lang-Trotter Conjectures on elliptic curves.

## Contents

## 0 Introduction

Let $A/L$ be an abelian variety of dimension $g$ over a number field. Deep arithmetic information about $A$ is encoded in its reductions $A_\mathfrak{p}$ modulo the primes $\mathfrak{p}$ of $L$. In particular, one may study the behavior of these reductions by tracking the properties of the the arithmetic Frobenius endomorphism $\pi_\mathfrak{p}$ of $A_\mathfrak{p}$ as $\mathfrak{p}$ varies with increasing norm.

The memoir of Lang-Trotter [LT76] studies the following questions for an elliptic curve $E/\mathbb{Q}$: how often does the Frobenius endomorphism of $E_p$ have a specified trace, $t \in \mathbb{Z}$? How often does $E_p$ have a specified imaginary quadratic number field $K$ as its endomorphism algebra? They give the following conjectures, which remain open:

**Conjecture 0.1** (Lang-Trotter Conjectures).

1. *Suppose that* $\operatorname{End}\left(E_{\overline{\mathbb{Q}}}\right) \cong \mathbb{Z}$, *or* $t \neq 0$. *Then there exists a constant* $C_{E,t} \geq 0$ *such that*

$$\#\left\{p \leq X \text{ of good reduction} \;\Big|\; \operatorname{Tr}(\pi_p) = t\right\} \sim C_{E,t} \frac{\sqrt{X}}{\log X},$$

*where* $C_{E,t} = 0$ *is understood to mean that the set written above is finite.*

---
* University of Maryland, College Park, MD, USA. *email:* bloom@math.umd.edu. *website:* www.math.umd.edu/~bloom

2. Suppose that $\mathrm{End}\left(E_{\overline{\mathbb{Q}}}\right) \cong \mathbb{Z}$, and let $K/\mathbb{Q}$ be an imaginary quadratic number field. Then there exists a constant $C_{E,K} > 0$ such that

$$\#\left\{p \leq X \text{ of good reduction} \;\middle|\; \mathrm{End}\left(E_p\right) \otimes \mathbb{Q} = \mathbb{Q}(\pi_p) \cong K\right\} \sim C_{E,K} \frac{\sqrt{X}}{\log X}.$$

The constants $C_{E,t}$ and $C_{E,K}$ have precise descriptions in terms of the statistical heuristics used and the Chebotarev density theorem. Many subsequent authors have studied variants of these questions, which we call under one umbrella "questions of Lang-Trotter type" (after [AH15]), namely:

**Question 0.2** (Question of Lang-Trotter Type)**.** *Let $A/L$ be an abelian variety of dimension $g$ over a number field. Let ♣ be a property of abelian varieties of dimension $g$ over finite fields. Give qualitative or numerical (i.e. asymptotic) descriptions of*

$$\Pi(A, \clubsuit) := \left\{\mathfrak{p} \subset \mathcal{O}_L \text{ of good reduction} : A_\mathfrak{p} \text{ has } \clubsuit\right\}.$$

By "qualitative descriptions" we mean descriptions of the primes which appear in $\Pi(A, \clubsuit)$ via congruence conditions, diophantine equations, and/or inequalities.

Our main result (Theorem 1.1) begins to answer a generalization of the second Lang-Trotter Conjecture. Namely, we study the question of Lang-Trotter type to describe those primes at which the reduction of a given abelian variety $A/\mathbb{Q}$ is ordinary, non-split, and has a specified CM field $K$ as its endomorphism algebra. Call the set of these primes $\Pi(A, K)$. With the assumption that the adelic Galois representation of $A$ has open image in $\mathrm{GSp}_{2g}\widehat{\mathbb{Z}}$, we find explicit asymptotic upper bounds

$$\#\left\{p \leq X \text{ of good reduction} \;\middle|\; p \in \Pi(A,K)\right\} \ll_{A,K} X(\log\log X)^\alpha/(\log X)^\beta, \qquad \text{unconditionally;}$$

$$\#\left\{p \leq X \text{ of good reduction} \;\middle|\; p \in \Pi(A,K)\right\} \ll_A X^\theta \log X, \qquad \text{under conjectural assumptions,}$$

where $\alpha, \beta > 1$, and $\theta < 1$, are functions of $g = \dim A$, and $\theta$ decreases with the strength of the conjectural assumptions. Our method, which is heavily inspired by [CFM05] and [Coj+16], uses a sieve (see Theorem 3.1) and relies on the properties of the Galois representations of $A$ to bound the size of $\Pi(A,K)$. Because of this, we must restrict to the proper class of varieties whose adelic Galois representation is "eventually surjective," as stated above. In particular, this hypothesis implies that $A$ is without extra endomorphisms.

## 0.1 Outline of the Article.

In the remainder of this Introduction, we introduce notations and more precisely introduce the question of Lang-Trotter type that we investigate. In Section 1, we present the main results of this paper, namely Theorem (1.1) and Theorem (1.2), as well as Corollaries of these results. In Section 2, we review the literature on questions of Lang-Trotter type. The reader may skip ahead after Section 1 to Section 3 without loss of continuity, although we refer in Section 6 to a few facts stated in Section 2. In Section 3, we review the preliminaries needed for our proofs. In Section 4, we prove Theorem 1.1. In Section 5, we prove Theorem 1.2. In Section 6, we prove the Corollaries stated in Section 1 In Section 7, we make concluding remarks and present questions for further study.

## 0.2 Notations.

We use the standard Bachmann-Landau notations for asymptotic growth of functions, which we now recall. A subscript $\star$ will denote that the implied constant depends *only* on the object(s) $\star$, so that if $\star$ is empty, then the implied constant is absolute. We write $X \gg_\star 0$ to mean "for all $X \geq N_\star$." Let $f, g : \mathbb{N} \to \mathbb{R}$. We write $g(X) = O_\star(f(X))$ or $g(X) \ll_\star f(X)$ to mean $\exists C_\star \geq 0$ such that for $X \gg_\star 0$, $|g(X)| \leq C_\star |f(X)|$. We write $g(X) = o(f(X))$ to mean $\lim_{X \to \infty} \frac{g(X)}{f(X)} = 0$, and we write $g(X) \sim f(X)$ to mean $\lim_{X \to \infty} \frac{g(X)}{f(X)} = 1$. We write $g(X) \asymp_\star f(X)$ to mean "$g(X) \ll_\star f(X)$ and $f(X) \ll_\star g(X)$."

We will use the letters $l$, $p$, $q$, and $\ell$ to denote rational prime numbers, and $\mathfrak{p}$ to denote a prime ideal in a number field. We will use N and Tr to denote "norm" and "trace", respectively, when the meaning is clear from context, and introduce subscripts and superscripts when the meaning is not clear. We will write $\left(\frac{\alpha}{a}\right)$ for the Jacobi (i.e., generalized quadratic residue) symbol of $\alpha$ modulo $a$. In a number field $L$, we will write $n_L$ for the degree of the extension $L/\mathbb{Q}$, $d_L$ for the discriminant of the extension $L/\mathbb{Q}$, and $h_L$ for the class number. For any set $S$ of prime ideals of a number field (or rational prime numbers), we denote the prime-counting function

$$S(X) := \# \left\{ \mathfrak{p} \in S \ \Big| \ \mathrm{N}\,\mathfrak{p} \leq X \right\}.$$

We say that $S$ has (natural) density $\delta$ if

$$\lim_{X \to \infty} \frac{S(X)}{\#\left\{\mathfrak{p} \ \Big| \ \mathrm{N}\,\mathfrak{p} \leq X\right\}} = \delta.$$

For a finite set $X$, we will write $\#X$ for the cardinality. For a finite group $G$ and a union of conjugacy classes $C \subseteq G$, we will write $\widetilde{C}$ for the number of conjugacy classes contained in $C$.

For an abelian variety over a field $A/\kappa$, we will always use $\mathrm{End}(A)$ to denote the ring of endomorphisms of $A$ defined over the base field $\kappa$.

### 0.3 Our Question.

As mentioned earlier, the question that we study here is an extension of the "fixed-field" Lang-Trotter question to higher-dimensional abelian varieties. Honda-Tate theory [Hon68; Tat69] tells us that when $p$ is a prime of good, ordinary, non-split reduction for $A$, then its endomorphism algebra $\mathrm{End}(A_p) \otimes \mathbb{Q}$ is a CM field of degree $2g$, equal to its **Frobenius field** $\mathbb{Q}(\pi_p)$. It is known as well that $\mathrm{End}(A)$ (the endomorphism ring from characteristic zero) embeds into $\mathbb{Q}(\pi_p)$. Thus, when $A$ does not have CM, its Frobenius fields are CM fields of degree $2g$ that admit an embedding of $\mathrm{End}(A)$ as a subring. We thus ask the following Question.

**Question 0.3.** *Let $A/\mathbb{Q}$ be a non-CM abelian variety of dimension $g$. Let $K$ be a CM field of degree $2g$. Describe*

$$\Pi(A, K) := \left\{ p \text{ of good, ordinary, nonsplit reduction} \ \Big| \ K \cong \mathrm{End}(A_p) \otimes \mathbb{Q} \right\}.$$

We also ask about supersets of $\Pi(A, K)$; namely, we ask

**Question 0.4.** *Let $A/\mathbb{Q}$ be a non-CM abelian variety of dimension $g$. Let $F$ be a totally real field of degree $g$. Describe*

$$\Pi(A, F) := \left\{ p \text{ of good, ordinary, nonsplit reduction} \ \Big| \ F \hookrightarrow \mathrm{End}(A_p) \otimes \mathbb{Q} \right\}.$$

### Acknowledgments

I would like to thank my advisor, Professor Larry Washington, for his advice and encouragement during this research. I would also like to thank Professor Washington and Professor Niranjan Ramachandran for reviewing earlier versions of this article.

## 1 Main Results

We mimic the application in [CFM05] of the Square Sieve (Theorem 3.1) to obtain the following.

**Theorem 1.1.** *Let $A/\mathbb{Q}$ be a principally polarized abelian variety of conductor $N$ whose adelic Galois representation $\widehat{\rho}$ has image that is open in $\mathrm{GSp}_{2g} \widehat{\mathbb{Z}}$. (See Section 3 for definitions and the Remark below.) Let $K/\mathbb{Q}$ be a CM field of degree $2g$ with discriminant $d = d(K/\mathbb{Q})$. Then,*

$$\Pi(A, K)(X) \ll_{N, g} \begin{cases} X^{1 - 1/(8g^2 + 4g + 6)} \log X & \text{under GRH;} \\ X^{1 - 1/(4g^2 + 4g + 6)} \log X & \text{under GRH and AHC;} \\ X^{1 - 1/(2g^2 + 4g + 6)} \log X & \text{under GRH, AHC, and PCC;} \\ X (\log \log X)^{1 + 1/(4g^2 + 3g + 2)} (\log X)^{-1 - 1/(8g^2 + 6g + 4)} (1 + \nu(d)) & \text{unconditionally,} \end{cases}$$

*where $\nu(d)$ is the number of distinct prime divisors of $d$. The conjectural assumptions are as follows:*

**GRH:** *the Generalized Riemann Hypothesis holds for the Dedekind zeta function of the division fields $\mathbb{Q}(A[lq])/\mathbb{Q}$, for all distinct primes $l, q \gg 0$;*

**AHC:** *Artin's Holomorphy Conjecture holds for the Artin L-functions attached to the irreducible characters of $\operatorname{Gal}\mathbb{Q}(A[lq])/\mathbb{Q}$, for all distinct primes $l, q \gg 0$;*

**PCC:** *a certain Pair Correlation Conjecture holds for the Artin L-functions attached to the irreducible characters of $\operatorname{Gal}\mathbb{Q}(A[lq])/\mathbb{Q}$, for all distinct primes $l, q \gg 0$.*

See [Mur01] for precise formulations of Conjectures AHC and PCC.

**Theorem 1.2.** *Let $A/\mathbb{Q}$ be a principally polarized abelian surface with $\operatorname{End}\left(A_{\overline{\mathbb{Q}}}\right) \cong \mathbb{Z}$. Let $F = \mathbb{Q}(\sqrt{d})$ be a real quadratic number field, where $d$ is squarefree. Then,*

$$\Pi(A, F)(X) \ll_N \begin{cases} X^{45/46} \log X & \text{under GRH;} \\ X^{29/30} \log X & \text{under GRH and AHC;} \\ X^{22/23} \log X & \text{under GRH, AHC, and PCC;} \\ X(\log \log X)^{23/22}(\log X)^{-67/66}\left(1 + \nu(d)\right) & \text{unconditionally,} \end{cases}$$

*where $\nu(d)$ is the number of distinct prime divisors of $d$. The conjectural assumptions are identical to those above.*

**Remark 1.3.** *The hypothesis in Theorem 1.1 that $\operatorname{im}\widehat{\rho}$ be open in $\operatorname{GSp}_{2g}\widehat{\mathbb{Z}}$ implies that $A$ without extra endomorphisms, i.e, $\operatorname{End}\left(A_{\overline{\mathbb{Q}}}\right) \cong \mathbb{Z}$. Moreover, the hypothesis is true for a wide class of varieties without extra endomorphisms. Works of Serre [Ser00b; Ser00a] and Pink [Pin98] show that the hypothesis is true when $\operatorname{End}\left(A_{\overline{\mathbb{Q}}}\right) \cong \mathbb{Z}$ if $g = 1, 2$ or if $g \geq 3$ is not in the set*

$$\left\{\frac{1}{2}(2n)^k \,\middle|\, n > 0, k \geq 3 \text{ odd}\right\} \cup \left\{\frac{1}{2}\binom{2n}{n} \,\middle|\, n \geq 3 \text{ odd}\right\} = \{4, 10, 16, 32, \ldots\}$$

*The hypothesis is also true for those p.p.a.v. satisfying the property "(T)" of [Hal11]. Thus, adding this hypothesis to Theorem 1.2 would be redundant.*

We also consider the set of those CM fields which appear as Frobenius fields of $A$,

$$\mathcal{D}_A := \left\{\mathbb{Q}(\pi_p) \,\middle|\, p \text{ good, ordinary, non-split}\right\}$$

and the set of their totally real subfields,

$$\mathcal{D}_A^0 := \left\{\mathbb{Q}(\pi_p)_0 \,\middle|\, p \text{ good, ordinary, non-split}\right\}.$$

If $A$ is a surface, we index (essentially) by discriminant,

$$\mathcal{D}_A^0(X) := \left\{\mathbb{Q}(\sqrt{d}) \in \mathcal{D}_A^0 \,\middle|\, d \text{ squarefree}, d \leq 48X\right\}$$

For an abelian variety $A$ of dimension $g$, we index $\mathcal{D}_A$ by certain effective functions $\psi_g(\sqrt{X})$ which are polynomials in $\sqrt{X}$:

$$\mathcal{D}_A(X) := \left\{K \in \mathcal{D}_A \,\middle|\, \operatorname{sf}(d(K/\mathbb{Q})) \leq \psi_g(\sqrt{X})\right\}$$

where $\operatorname{sf}(d)$ is the square-free part of $d$. See the discussion after Corollary 4.5 for details.

Using the Pigeonhole Principle, we obtain from our main Theorems the following asymptotic lower bounds on the size of $\#\mathcal{D}_A(X)$ and $\#\mathcal{D}_A^0(X)$.

**Corollary 1.4.** *Let the notations be as above. Let $\delta$ be the density of the set of good, ordinary, non-split primes for $A$. If $g > 2$, assume that $\delta > 0$. Then,*

$$\#\mathcal{D}_A(X) \gg_N \delta \frac{X^\theta}{(\log X)^2},$$

*where we may take*

$$\theta = \begin{cases} 1/(8g^2 + 4g + 6) & \text{under GRH;} \\ 1/(4g^2 + 4g + 6) & \text{under GRH and AHC;} \\ 1/(2g^2 + 4g + 6) & \text{under GRH, AHC, and PCC.} \end{cases}$$

**Corollary 1.5.** *Let the notations be as above, and suppose that $A$ is a surface. Then,*

$$\#\mathcal{D}_A^0(X) \gg_N \frac{X^\theta}{(\log X)^2},$$

*where we may take*

$$\theta = \begin{cases} 1/46 & \text{under GRH;} \\ 1/30 & \text{under GRH and AHC;} \\ 1/23 & \text{under GRH, AHC, and PCC.} \end{cases}$$

**Corollary 1.6.** *Let the notations be as above. Unconditionally, $\#\mathcal{D}_A(X) \to \infty$, and if $A$ is a surface, $\#\mathcal{D}_A^0(X) \to \infty$.*

## 2 Background

For the Subsections below, we let $A/L$ be an absolutely simple abelian variety *without Complex Multiplication (non-CM)* of dimension $g$ over a number field; that is to say, $A_{\overline{L}}$ is not isogenous to a product of abelian varieties of smaller dimension, and $\text{End}(A_{\overline{L}}) \otimes \mathbb{Q}$ is *not* a number field of degree $2g$. We also let $E/L$ be a *non-CM* elliptic curve over a number field. In either context, $N$ is the conductor. We also let $B/\mathbb{F}_p$ be an abelian variety of dimension $g$.

### 2.1  $p$-Rank.

Recall that the group of geometric $p$-torsion of $B$ has the shape

$$B(\overline{\mathbb{F}}_p) \cong (\mathbb{Z}/p\mathbb{Z})^f$$

for some $0 \leq f \leq g$. We call the integer $f$ the **p-rank** of $B$. If $f = g$, we call $B$ **ordinary**, otherwise we call $B$ **non-ordinary**. If $B$ is an elliptic curve or abelian surface with $f = 0$, we call $B$ **supersingular**.[1]

It is known that, possibly only after a finite extension of the base-field $L$ of $A$, the set of non-ordinary primes $\Pi(A, f \neq g)$ has density zero if $g = 1$ [Ser68], if $g = 2$ [Ogu81], and for some abelian varieties with $g = 3$ [Tan99] or $g$ a power of 4 [Noo95]. For arbitrary $g$, $\Pi(A, f \geq 2)$ has density one [Ogu81; BG97], but it is not known in general whether the set of ordinary primes for $A$ has positive density. Because $E_p$ is supersingular iff its trace of Frobenius is 0 (for $p \geq 5$), the Conjecture 0.1 predicts the asymptotics of $\Pi(E, f = 0)$ for $E/\mathbb{Q}$.

Various authors have improved upon the upper bound for $E/L$ of [Ser68]. The best known upper bounds for $E/\mathbb{Q}$ are

$$\Pi(E, f = 0)(X) \ll_N \begin{cases} X^{3/4} & \text{unconditionally [Elk87b; Elk87a];} \\ X^{3/4}(\log X)^{-1/2} & \text{under GRH [Zyw15]} . \end{cases}$$

---

[1] This adjective means that the Newton slopes at $p$ of the characteristic polynomial of $\pi_p$ are all $1/2$. This condition is equivalent to $f = 0$ only when $g \leq 2$.

(It is astounding how small an improvement GRH affords with our current technology!) As for lower bounds, [Elk87c; Elk89] prove that if $L$ has a real embedding (e.g., if $L = \mathbb{Q}$), $\Pi(E, f = 0)$ is infinite. For $L = \mathbb{Q}$, various authors improve this lower bound; the best known bound is

$$\Pi(E, f = 0)(X) \gg_N \begin{cases} \log \log X & \text{under GRH [Elk87b];} \\ \frac{\log \log \log X}{(\log \log \log \log X)^{1+\epsilon}} & \text{unconditionally [FM96].} \end{cases}$$

Much less is known about higher-dimensional non-CM abelian varieties $A$. The author knows of no bounds better than $\Pi(A, f \neq g) = o(\pi(X))$ for only those abelian varieties mentioned in the second paragraph of this Subsection, and he knows of no asymptotic lower bounds if $g \geq 2$, even for a single non-CM abelian variety. Nor is it known whether $\#\Pi(A, f \neq g) = \infty$ for any non-CM abelian variety of dimension $g \geq 2$.

If $A$ **has real multiplication**, which means here that $\text{End}(A) \otimes \mathbb{Q}$ is a totally real number field of degree $g$ and $\text{End}(A) \otimes \mathbb{Q} \cong \text{End}\left(A_{\overline{\mathbb{Q}}}\right) \otimes \mathbb{Q}$, [BG97] conjectures a probabilistic model which yields

$$\Pi(A, f < g)(X) \sim \begin{cases} C_A \frac{\sqrt{X}}{\log X} & \text{if } g = 1, \\ C_A \log \log X & \text{otherwise,} \end{cases}$$

and

$$\Pi(A, f = 0)(X) \sim \begin{cases} C_{A,0} \frac{\sqrt{X}}{\log X} & \text{if } g = 1, \\ C_{A,0} \log \log X & \text{if } g = 2, \\ O(1) & \text{otherwise,} \end{cases}$$

for certain positive constants $C_A$ and $C_{A,0}$. This conjecture remains open.

## 2.2 Fixed-Trace.

For an elliptic curve $E/\mathbb{Q}$, denote $a_p := \text{Tr}(\pi_p)$. For primes $p$ of good reduction, i.e. if $E_p$ is an elliptic curve, then $\#E_p(\mathbb{F}_p) = p + 1 - a_p$. The Hasse-Weil bound states that $|a_p| \leq 2\sqrt{p}$.

For a fixed integer $t \neq 0$, Conjecture 0.1 predicts the size of $\Pi(E, a_p = t)(X)$. Various upper and lower bounds (conditional and unconditional) are given in the literature for $\Pi(E, a_p = t)(X)$. (We restrict to $t \neq 0$ because $E_p$ is supersingular iff $a_p = 0$ when $p \geq 5$.) Unconditionally, Serre [Ser81] gives the first bound, $\Pi(E, a_p = t)(X) \ll_N X/(\log X)^{5/4-\epsilon}$. This was improved by Wan [Wan90] and Murty [Mur97]. The best known unconditional upper bound is from the recent preprint [TZ16], which gives

$$\Pi(E, a_p = t)(X) \ll_N \frac{X(\log \log X)^2}{(\log X)^2}.$$

Conditionally on GRH, Serre [Ser81] also gives the first bound, $\Pi(E, a_p = t)(X) \ll_N X^{7/8}(\log X)^{1/2}$. This was improved by Murty-Murty-Saradha [MMS88]. The best known upper bound (conditional on GRH) is

$$\Pi(E, a_p = t)(X) \ll_N X^{4/5}(\log X)^{-3/5}$$

of Zywina [Zyw15].

For higher-dimensional abelian varieties, this question has just begun investigation. The recent work of Cojocaru-Davis-Silverberg-Stange [Coj+16] studies the $\text{GL}_{2g}$-trace of Frobenius, $a_{1,p} := \text{Tr}\,\pi_p$ for the class of abelian varieties $A/\mathbb{Q}$ whose adelic Galois representation $\widehat{\rho}$ (see (2)) has open image in $\text{GSp}_{2g}\widehat{\mathbb{Z}}$. They obtain the bounds

$$\Pi(A, a_{1,p} = t)(X) \ll_{A,\epsilon} \begin{cases} X^{1-\frac{1}{2}\theta+\epsilon} & \text{under GRH,} \\ X/(\log X)^{1+\theta-\epsilon} & \text{unconditionally,} \end{cases}$$

where $0 < \theta < 1/4$ decreases as $g$ increases. They obtain improvements upon the above, in the form of larger $\theta$, when $t \neq 2g$, and further improvements when $t = 0$. Moreover, they show that with a conjectural assumption on the behavior of the Galois representations of $A$ that generalizes the Sato-Tate Conjecture,

$$\Pi(A, a_{1,p} = t)(X) \sim C_{A,t}\frac{\sqrt{X}}{\log X}$$

for some precisely defined constant $C_{A,t} \geq 0$, where, as before, we understand $C_{A,t} = 0$ to mean that the set is finite.

## 2.3 Fixed-Field.

For any elliptic curve $E/\mathbb{Q}$ and a prime $p$ of good reduction, it is well-known that the endomorphism algebra $\operatorname{End}(E_p) \otimes \mathbb{Q} = \mathbb{Q}(\pi_p) \cong \mathbb{Q}(\sqrt{D_p})$ is an imaginary quadratic field, where we take $D_p$ to be the squarefree part of $a_p^2 - 4p$. If $E_p$ is supersingular, so that $a_p = 0$ (for $p \geq 5$), then $\operatorname{End}(E_p) \cong \mathbb{Q}(\sqrt{-p})$, and $\operatorname{End}\left((E_p)_{\overline{\mathbb{F}}_p}\right)$ is isomorphic to the quaternion algebra $\mathbb{B}_{p,\infty}/\mathbb{Q}$ ramified only at $p$ and $\infty$. But if $E_p$ is ordinary, then $D_p$ might possibly take any squarefree value between $0$ and $-4p$, not inclusive, and $E_p$ does not pick up any extra endomorphisms over $\overline{\mathbb{F}}_p$.

In particular, if $E$ is *non-CM*, then the endomorphism algebras (or **Frobenius fields**) at ordinary primes vary in the set of imaginary quadratic number fields $K$. The article of Cojocaru-Fouvry-Murty [CFM05] investigates the sets

$$\Pi(E, K) \coloneqq \Pi\left(E, \ p \text{ ordinary and } \operatorname{End}(E_p) \otimes \mathbb{Q} \cong K\right)$$

via the Square Sieve (see Subsection 3.1) and obtains the first bounds in print. These bounds are of the form $\Pi(E,K)(X) \ll_N X^\theta \log X$, conditional on various conjectural assumptions, and $\ll_{N,d(K/\mathbb{Q})} (\log \log X)^{13/12} (\log X)^{-25/24}$ unconditionally. See the remarks preceding the statement of this Theorem in [CFM05] for a history of remarks made by other authors which indicated bounds on $\Pi(E,K)(X)$. Improvements on these bounds have been made by various authors [CD08; Zyw15; TZ16] using sieves and "mixed representations" as suggested by Serre. The best known upper bounds are

$$\Pi(E,K)(X) \ll_E X^{4/5}(\log X)^{-3/5} h_K^{-3/5} + X^{1/2}(\log X)^3, \qquad \text{under GRH [Zyw15];}$$
$$\Pi(E,K)(X) \ll_{E,K} X(\log \log X)(\log X)^{-2}, \qquad \text{unconditionally [TZ16].}$$

## 2.4 Geometrically Simple.

Suppose $A$ is geometrically simple, i.e. $A_{\overline{L}}$ is simple. Murty-Patankar [MP08] investigated the set of primes $\Pi(A, \text{geom. simple})$ at which $A$ remains geometrically simple. They show that if $A$ has Complex Multiplication or has Real Multiplication, then $\Pi(A, \text{geom. simple})$ has density one. Moreover, they and Zywina [Zyw13] conjecture that for any $A/L$, possibly[2] after a finite extension $L'/L$,

$$\operatorname{End}(A_{\overline{L}}) \text{ is commutative} \iff \Pi(A_{L'}, \text{geom. simple}) \text{ has density } \delta_{A_{L'}} = 1 \tag{1}$$

Achter [Ach09] proves the backward direction of (1), and shows that moreover, if $\operatorname{End}(A_{\overline{L}})$ is non-commutative, there is a finite extension $L'/L$ such that $\delta_{A_{L'}} = 0$. He moreover also proves the forward direction of (1) if $\operatorname{End}(A_{\overline{L}})$ is a totally real or totally imaginary field and if $A$ satisfies a certain parity assumption. The particular case of the forward direction of (1) when $\operatorname{End}(A_{\overline{L}}) \cong \mathbb{Z}$ is an earlier result of Chavdarov [Cha97]. Achter [Ach12] gives explicit bounds on $\Pi(A, \text{geom. split})(X)$ in these cases (and one other). Zywina [Zyw13] proves that if the Manin-Mumford conjecture is true for $A$, then possibly after a finite extension, the forward direction is true. Murty-Zong [MZ14] prove that if for some prime $\ell \in \mathbb{Z}$, $\operatorname{End}(A) \otimes \mathbb{Q}_\ell \cong \operatorname{End}(A_{\overline{L}}) \otimes \mathbb{Q}_\ell$ is a field, if the Zariski closure of the image of the $\ell$-adic Galois representation $\rho_{\ell^\infty}$ is connected, and if $\rho_{\ell^\infty}$ satisfies an additional technical assumption, then $\delta_A > 0$. Lastly, we mention the preprint [AH15] which estimates the number of *split* abelian surfaces over $\mathbb{F}_p$ as approximately $p^{-1/2}\left(\#\mathcal{A}_2 \mathbb{F}_p\right)$, from which they conjecture that for an abelian surface $A/\mathbb{Q}$ without extra endomorphisms, $\Pi(A, \text{split})(X) \sim C_A \frac{\sqrt{X}}{\log X}$ for some positive constant $C_A$.

## 3 Preliminaries

### 3.1 The Square Sieve

As in [CFM05], the main tool we use is the square sieve, which originates in [Hea84].

**Theorem 3.1** (Square Sieve). *Let $\mathcal{A}$ be a finite sequence of non-zero rational integers, and $\mathcal{P}$ a set of distinct odd rational primes. Set*

$$S(\mathcal{A}) \coloneqq \#\{\alpha \in \mathcal{A} : \alpha \text{ is a square}\}.$$

---

[2] As pointed out in [Zyw13], there are counterexamples to the conjecture without the extension.

*Then,*

$$S(\mathcal{A}) \leq \frac{\#\mathcal{A}}{\#\mathcal{P}} + \max_{\substack{l,q \in \mathcal{P} \\ l \neq q}} \left| \sum_{\alpha \in \mathcal{A}} \left(\frac{\alpha}{lq}\right) \right| + \frac{2}{\#\mathcal{P}} \sum_{\alpha \in \mathcal{A}} \sum_{\substack{l \in \mathcal{P} \\ (\alpha, l) \neq 1}} 1 + \frac{1}{(\#\mathcal{P})^2} \sum_{\alpha \in \mathcal{A}} \left( \sum_{\substack{l \in \mathcal{P} \\ (\alpha, l) \neq 1}} 1 \right)^2$$

*where $\left(\frac{\cdot}{\cdot}\right)$ is the Jacobi symbol.*

*Proof.* See, for instance, Section 2.1 of [CFM05]. □

## 3.2 Explicit Chebotarev Density Theorems

Let $L/\mathbb{Q}$ be a finite Galois extension with Galois group $G$, degree $n_L$, and discriminant $d_L$. Let $C$ be a union of conjugacy classes of $G$. Denote by $\mathcal{P}(L/\mathbb{Q})$ the set of rational primes $p$ which ramify in $L/\mathbb{Q}$. Set

$$M(L/\mathbb{Q}) := (\#G) \prod_{p \in \mathcal{P}(L/\mathbb{Q})} p.$$

Define the prime counting function for $C$,

$$\pi_C(X, L/\mathbb{Q}) := \# \left\{ p \leq X : p \text{ unramified in } L/\mathbb{Q}; \sigma_p \subseteq C \right\}$$

where $\sigma_p := \left(\frac{L/\mathbb{Q}}{p}\right)$ is the Artin symbol of $p$ in $L/\mathbb{Q}$. Recall that the Chebotarev density theorem states that as $x \to \infty$,

$$\pi_C(X, L/\mathbb{Q}) \sim \frac{\#C}{\#G} \int_2^X \frac{dt}{\log t}.$$

We use the notation $\operatorname{li} X := \int_2^X \frac{dt}{\log t}$ for the logarithmic integral to $X$. We will use "explicit" versions of this theorem; that is, versions with bounds on the error term of the approximation.

**Theorem 3.2** ([LO77; Ser81; MMS88; Mur97]). *Let the notation be as above. Then, for $X \gg 0$,*

$$\pi_C(X, L/\mathbb{Q}) = \frac{\#C}{\#G} \operatorname{li} X + R_C(X)$$

*where the error term $R_C(X)$ satisfies the following bounds:*

1. *Assume GRH for the Dedekind zeta function of $L/\mathbb{Q}$. Then,*

$$R_C(X) = O\left( (\#C) X^{1/2} \left( \frac{\log|d_L|}{n_L} + \log X \right) \right)$$

2. *Assume GRH and AHC for $L/\mathbb{Q}$. Then,*

$$R_C(X) = O\left( (\#C)^{1/2} X^{1/2} \left( \log M(L/\mathbb{Q}) + \log X \right) \right)$$

3. *Assume GRH, AHC, and PCC for $L/\mathbb{Q}$. Then,*

$$R_C(X) = O\left( (\#C)^{1/2} X^{1/2} \left( \frac{\#\widetilde{G}}{\#G} \right)^{1/4} \left( \log M(L/\mathbb{Q}) + \log X \right) \right)$$

4. *Unconditionally, there exist positive constants $A, B, B'$ with $A$ effective and $B, B'$ absolute, such that if*

$$\log X \geq B'(\#G) \left( \log|d_L| \right)^2,$$

   *then*

$$R_C(X) \ll \frac{\#C}{\#G} \operatorname{li}\left( X \exp\left( -B \frac{\log X}{\max\{|d_L|^{1/n_L}, \log|d_L|\}} \right) \right) + (\#\widetilde{C}) X \exp\left( -A \sqrt{\frac{\log X}{n_L}} \right),$$

*In all of the above, the implied constants are absolute.*

We will also employ the following bound on $|d_L|$ from [Ser81].

**Lemma 3.3.** *Let the notation be as above. Then,*

$$\frac{n_L}{2} \sum_{p \in \mathcal{P}(L/\mathbb{Q})} \log p \leq \log|d_L| \leq (n_L - 1) \sum_{p \in \mathcal{P}(L/\mathbb{Q})} \log p + n_L \log n_L.$$

## 3.3 Galois Representations and Open Image Varieties

Let $G_{\mathbb{Q}} := \mathrm{Gal}\left(\overline{\mathbb{Q}}/\mathbb{Q}\right)$. Let $A/\mathbb{Q}$ be a principally polarized abelian variety ("p.p.a.v.") of dimension $g$. Recall that for any integer $M \geq 1$, the geometric torsion subgroup

$$A[M](\overline{\mathbb{Q}}) \cong \left(\mathbb{Z}/M\mathbb{Z}\right)^{2g}$$

is naturally a $G_{\mathbb{Q}}$-module by action on the coordinates,

$$\rho_M : G_{\mathbb{Q}} \to \mathrm{GL}\left(A[M](\overline{\mathbb{Q}})\right) \cong \mathrm{GL}_{2g}\left(\mathbb{Z}/M\mathbb{Z}\right),$$

after choosing a basis of $A[M](\overline{\mathbb{Q}})$. However, the Galois action respects the Weil pairing $e_M$ on $A[M]$, so that in fact

$$\rho_M : G_{\mathbb{Q}} \to \mathrm{GSp}\left(A[M](\overline{\mathbb{Q}}), e_M\right) \cong \mathrm{GSp}_{2g}\left(\mathbb{Z}/M\mathbb{Z}\right),$$

after choosing a symplectic basis with respect to the Weil pairing. We call $\rho_M$ the **mod-$M$ Galois representation** of $A$. Let $\ell$ be a rational prime. We define the **$\ell$-adic Galois representation** as the inverse limit

$$\rho_{\ell^\infty} := \varprojlim \rho_{\ell^n} : G_{\mathbb{Q}} \to \mathrm{GSp}_{2g} \mathbb{Z}_\ell$$

and the **adelic Galois representation**

$$\widehat{\rho} := \prod_\ell \rho_{\ell^\infty} : G_{\mathbb{Q}} \to \prod_\ell \mathrm{GSp}_{2g} \mathbb{Z}_\ell \cong \mathrm{GSp}_{2g} \widehat{\mathbb{Z}} \tag{2}$$

The representations $\rho_M$, $\rho_{\ell^\infty}$, and $\widehat{\rho}$ are extremely important objects in the study of $A$. It is well-known that for $p \nmid N$ fixed and $\ell \neq p$ varying, the characteristic polynomial of Frobenius, $\mathrm{char}\, \rho_\ell(\pi_p) \in \mathbb{Z}[x]$, is independent of $\ell$. We will thus without comment use the notation $\mathrm{char}\, \pi_p$ or $\mathrm{char}_p$ for $\mathrm{char}\, \rho_\ell(\pi_p)$.

As stated in the Introduction, we study here those p.p.a.v. whose adelic representation $\widehat{\rho}$ has open image in $\mathrm{GSp}_{2g} \widehat{\mathbb{Z}}$. That is, we study $A$ such that for $\ell \gg_A 0$,

$$\mathrm{im}\, \rho_{\ell^\infty} \cong \mathrm{GSp}_{2g} \mathbb{Z}_\ell. \tag{3}$$

For the curiosity of the reader, we mention that it is a very hard open problem to remove the dependency on $A$ in the quantifier "$\ell \gg_A 0$" of the "open-image" results mentioned in Remark 1.3. That is to say, it is not currently known whether there is a uniform bound $\ell \gg_g 0$ such that (3) (or an appropriate modification thereof) holds for *every* p.p.a.v. of dimension $g$. This problem is known as the **Serre uniformity conjecture**. We also mention [Lom15a] and the recent preprint [Lom15b] which give explicit bounds, in terms of $g$ and the stable Faltings height of $A$, on the quantifier "$\ell \gg_A 0$" of these results.

## 3.4 The Lang-Weil Bound

We include here the bound of Lang-Weil [LW54] on the number of rational points of a variety over a finite field. We will employ this bound in the proof of our main Theorems.

**Theorem 3.4** ([LW54]). *Let $V \hookrightarrow \mathbb{P}^n_{\mathbb{F}_q}$ be a projective variety of dimension $r$ and degree $d$ over a finite field. Then,*

$$\left|\#V(\mathbb{F}_q) - q^r\right| = (d-1)(d-2)q^{r-\frac{1}{2}} + O_{n,r,d}(q^{r-1}).$$

We note in passing that this nearly recovers the Weil bound for the number of points on an abelian variety over a finite field.

## 3.5 Bounds on the size of sets in $\mathrm{GSp}_{2g}$.

In the proof of our main Theorems, we will employ a bound on the size of particular subsets of $\mathrm{GSp}_{2g} \mathbb{Z}/l\mathbb{Z}$. The bound appears (essentially) as stated below in [AH03] and originates in [Cha97].

We first recall a few well-known facts. For a prime $l$,

$$\# \mathrm{Sp}_{2g} \mathbb{F}_l = l^{g^2} \prod_{i=1}^{g} \left( l^{2i} - 1 \right) = l^{2g^2+g} - l^{2g^2+g-2} + O_g\left( l^{2g^2+g-6} \right) \tag{4}$$

There is the exact sequence $1 \to \mathrm{Sp}_{2g} \mathbb{F}_l \to \mathrm{GSp}_{2g} \mathbb{F}_l \xrightarrow{\mu} \mathbb{G}_m \mathbb{F}_l \to 1$, where $\mu$ is the **multiplicator character**, namely,

$$MJM^t = \mu(M)J$$

where $J = \begin{pmatrix} 0 & I_g \\ -I_g & 0 \end{pmatrix}$ is the matrix for the standard symplectic form. Thus,

$$\# \mathrm{GSp}_{2g} \mathbb{F}_l = (l-1)l^{g^2} \prod_{i=1}^{g} \left( l^{2i} - 1 \right) = l^{2g^2+g+1} - l^{2g^2+g} + O_g\left( l^{2g^2+g-1} \right) \tag{5}$$

Now let $f \in \mathbb{F}_l[x]$ be a characteristic polynomial of some matrix in $\mathrm{GSp}_{2g} \mathbb{F}_l$, and let $\mathrm{char}(M)$ denote the characteristic polynomial of $M$. Let

$$C(\mathbb{F}_l) := \left\{ M \in \mathrm{GSp}_{2g} \mathbb{F}_l \,\middle|\, \mathrm{char}\, M = f \right\}$$

be the set of matrices with specified characteristic polynomial $f$. ($C(\mathbb{F}_l)$ is the set of $\mathbb{F}_l$-valued points of a subscheme of $\mathrm{GSp}_{2g}/\mathbb{F}_l$, hence the notation.) Then,

**Lemma 3.5** ([Cha97]).

$$\frac{l^{2g^2}}{(l-1)(l+1)^{2g^2+g}} \leq \frac{\#C(\mathbb{F}_l)}{\# \mathrm{GSp}_{2g} \mathbb{F}_l} \leq \frac{l^{2g^2}}{(l-1)(l-1)^{2g^2+g}}. \tag{6}$$

This immediately implies that $\#C(\mathbb{F}_l) \asymp_g l^{2g^2}$. A form in which this Lemma will be useful to us is to consider the error term

$$Q_C := \frac{\#C(\mathbb{F}_l)}{\# \mathrm{GSp}_{2g} \mathbb{F}_l} - \frac{l^{2g^2}}{(l-1)(l+1)^{2g^2+g}}$$

which, by the above, satisfies

$$0 \leq Q_C \leq \frac{l^{2g^2}}{(l-1)(l-1)^{2g^2+g}} - \frac{l^{2g^2}}{(l-1)(l+1)^{2g^2+g}} \tag{7}$$

$$= \frac{l^{2g^2}}{(l-1)} \cdot \frac{2l\left( (l+1)^{(2g+1)(g-1)} + \cdots + (l-1)^{(2g+1)(g-1)} \right)}{(l^2-1)^{2g^2+g}} \tag{8}$$

$$\ll l^{2g^2-1+\left(1+(2g+1)(g-1)\right)-2(2g^2+g)} = l^{-3g-1}. \tag{9}$$

We will also need to bound the number of conjugacy classes in $\mathrm{GSp}_{2g} \mathbb{Z}/lq\mathbb{Z}$, i.e. $\#\widetilde{\mathrm{GSp}_{2g} \mathbb{Z}/lq\mathbb{Z}}$. The paper [FG12], based on work of Wall [Wal63], gives the following bounds.

**Lemma 3.6.** *Let $g \geq 1$. Then, $q^g \leq \#\widetilde{\mathrm{Sp}_{2g} \mathbb{Z}/q\mathbb{Z}} \leq 10.8 q^g$*

With the Chinese Remainder Theorem and the long exact sequence of Section (1.4) of [HK85] that relates $\widetilde{\mathrm{GSp}_{2g} \mathbb{Z}/lq\mathbb{Z}}$ with $\widetilde{\mathrm{Sp}_{2g} \mathbb{Z}/lq\mathbb{Z}}$ and $\widetilde{\mathbb{F}_l^\times}$, we may thus conclude that

**Lemma 3.7.** $\#\widetilde{\mathrm{GSp}_{2g} \mathbb{Z}/lq\mathbb{Z}} \ll l^{g+1} q^{g+1}$.

# 4 Proof of Theorem 1.1

In this section, $A/\mathbb{Q}$ is a principally polarized abelian variety of dimension $g$ whose adelic Galois representation $\hat{\rho}$ has open image in $\mathrm{GSp}_{2g}\hat{\mathbb{Z}}$. This implies that $A$ is simple. Let $N$ be the conductor of $A$.

Let $p \nmid N$ be a prime of good, ordinary, non-split reduction for $A$. Then, by Honda-Tate theory, the endomorphism algebra $K := \mathrm{End}(A_p) \otimes \mathbb{Q}$ is a CM field of degree $2g$. Let $K_0$ be the totally real subfield of $K$. Then,

$$K \cong \mathbb{Q}(\pi_p) = K_0(\sqrt{r})$$

for some totally negative integer $r \in K_0$.

Because $\pi_p$ is a $p$-Weil number, the characteristic polynomial of the Frobenius endomorphism $\pi_p$ has the shape

$$\mathrm{char}_p(x) = x^{2g} + a_{1,p} x^{2g-1} + \ldots + a_{g,p} x^g + p a_{g-1,p} x^{g-1} + \ldots + p^g, \tag{10}$$

and the Triangle Inequality yields

$$|a_{i,p}| \leq \binom{2g}{i} p^{i/2}.$$

For convenience, when the prime $p$ is clear from context, we suppress it from the subscripts.

The following Lemmas, specifically Corollary 4.5, allow us to apply the Square Sieve.

**Lemma 4.1.** *Let the notation be as above, but with $K$ an arbitrary CM field of degree $2g$. Then,*

$$\mathbb{Q}(\pi_p) \cong K \implies \mathrm{N}_{\mathbb{Q}}^{K_0}\left((\pi_p + \overline{\pi_p})^2 - 4p\right) \cdot d(K/\mathbb{Q}) \in \mathbb{Z}^2$$

*Proof.* Let $x \mapsto \overline{x}$ be the complex conjugation of $K/K_0$. Then,

$$K_0 = \mathbb{Q}(\pi + \overline{\pi}); \quad K = K_0(\pi)$$

so that the ideal $d(K/K_0)$ is equal (up to the square of an ideal) to the discriminant of $x^2 - (\pi + \overline{\pi})x + p$. That is, $d(K/K_0) \cdot \mathfrak{a}^2 = ((\pi_p + \overline{\pi_p})^2 - 4p)\mathcal{O}_{K_0}$ for some ideal $\mathfrak{a}$ of $K_0$. Then, the formula for the norm of the relative discriminant gives

$$d(K/\mathbb{Q}) = \mathrm{N}_{\mathbb{Q}}^{K_0}\left(\frac{((\pi_p + \overline{\pi_p})^2 - 4p)\mathcal{O}_{K_0}}{\mathfrak{a}^2}\right) d(K_0/\mathbb{Q})^{[K:K_0]}$$

so that

$$\mathrm{N}_{\mathbb{Q}}^{K_0}\left((\pi_p + \overline{\pi_p})^2 - 4p\right) = \frac{d(K/\mathbb{Q}) \, \mathrm{N}\,\mathfrak{a}^2}{d(K_0/\mathbb{Q})^2},$$

and the relation follows. □

We note in passing that the sign of both sides is $(-1)^g$, so that in fact the above is an equality in $\mathbb{Z}$ and not just of ideals.

**Lemma 4.2.** *Suppose the integer $\pi + \overline{\pi}$ has minimal polynomial over $\mathbb{Q}$ equal to $x^g + \sum_{j=0}^{g-1} c'_j x^{g-j}$. Then, the $c'_j$ are polynomials of the $a_i$ and of $p$. These polynomials depend only on $g$.*

*Proof.* If $x^g + \sum_{j=0}^{g-1} c'_j x^{g-j} = 0$ is the minimal polynomial of $\pi + \overline{\pi}$, then

$$0 = \left(\pi + \frac{p}{\pi}\right)^g + c'_1 \left(\pi + \frac{p}{\pi}\right)^{g-1} + \ldots + c'_g$$

$$= \sum_{j=0}^{g} c'_j \sum_{k=0}^{g-j} \binom{g-j}{k} \pi^k \left(\frac{p}{\pi}\right)^{g-j-k}$$

so, multiplying through by $\pi^g$,

$$0 = \sum_{j=0}^{g} c'_j \sum_{k=0}^{g-j} \binom{g-j}{k} \pi^{2k+j} p^{g-j-k} \tag{11}$$

The result follows from solving the system of equations that results from comparing (11) to the minimal polynomial of $\pi$. □

**Lemma 4.3.** Let $n \geq 1$. Let $A \in \mathrm{GL}_n \mathbb{Z}_\ell$ have characteristic polynomial $x^n + \sum_{i=0}^{n-1} \alpha_i x^{n-i}$. Then, the coefficients of the characteristic polynomial of $A^2$ are polynomials in the $\alpha_i$. These polynomials do not depend on $A$, and are at worst quadratic in each of the $\alpha_i$.

*Proof.* $A$ is similar to a matrix in "companion form,"

$$A \sim \begin{pmatrix} 0 & 0 & 0 & \cdots & 0 & -\alpha_0 \\ 1 & 0 & 0 & \cdots & 0 & -\alpha_1 \\ 0 & 1 & 0 & \cdots & 0 & -\alpha_2 \\ 0 & 0 & 1 & \cdots & 0 & -\alpha_3 \\ \vdots & \vdots & \vdots & \ddots & \vdots & \vdots \\ 0 & 0 & 0 & \cdots & 1 & -\alpha_{n-1} \end{pmatrix} \implies A^2 - xI \sim \begin{pmatrix} -x & 0 & 0 & \cdots & 0 & -\alpha_0 & \alpha_{n-1}\alpha_0 \\ 0 & -x & 0 & \cdots & 0 & -\alpha_1 & \alpha_{n-1}\alpha_1 - \alpha_0 \\ 1 & 0 & -x & \cdots & 0 & -\alpha_2 & \alpha_{n-1}\alpha_2 - \alpha_1 \\ 0 & 1 & 0 & \cdots & 0 & -\alpha_3 & \alpha_{n-1}\alpha_3 - \alpha_2 \\ \vdots & \vdots & \vdots & \ddots & \vdots & \vdots & \vdots \\ 0 & 0 & 0 & \cdots & 1 & -\alpha_{n-1} & \alpha_{n-1}^2 - \alpha_{n-2} - x \end{pmatrix}.$$

Perform the column operation adding $\alpha_n \cdot$ (column $n-1$) to column $n$:

$$\det(A^2 - xI) = \det \begin{pmatrix} -x & 0 & 0 & \cdots & 0 & -\alpha_0 & 0 \\ 0 & -x & 0 & \cdots & 0 & -\alpha_1 & -\alpha_0 \\ 1 & 0 & -x & \cdots & 0 & -\alpha_2 & -\alpha_1 \\ \vdots & \vdots & \vdots & \ddots & \vdots & \vdots & \vdots \\ 0 & 0 & 0 & \cdots & -x & -\alpha_{n-3} & -\alpha_{n-4} \\ 0 & 0 & 0 & \cdots & 0 & -\alpha_{n-2} - x & -\alpha_{n-3} - \alpha_{n-1}x \\ 0 & 0 & 0 & \cdots & 1 & -\alpha_{n-1} & -\alpha_{n-2} - x \end{pmatrix}$$

Expanding out the determinant of the right-hand side, we see that each term in $\det(A^2 - xI)$ is at worst quadratic in each $\alpha_i$. $\square$

**Lemma 4.4.** Suppose the integer $\beta := (\pi + \overline{\pi})^2$ has characteristic equation

$$\mathrm{charpoly}_\beta(x) := x^g + c_1 x^{g-1} + \ldots + c_g := \prod_{\tau: K_0 \hookrightarrow \overline{\mathbb{Q}}} (x - \tau(\beta)) = 0$$

when considered as a linear transformation on the vector space $\mathbb{Q}(\pi + \overline{\pi})$ over $\mathbb{Q}$. (That is, considered as the multiplication map $x \mapsto \beta x$.) Then, the $c_i$ are polynomials of the $a_i$ and $p$, and these polynomials depend only on $g$. Moreover, these polynomials are at worst quadratic in the $a_i$.

*Proof.* This follows from the previous two Lemmas. $\square$

Thus, by noting that

$$\mathrm{N}_\mathbb{Q}^{K_0}\left((\pi_p + \overline{\pi_p})^2 - 4p\right) = (-1)^g \cdot \mathrm{charpoly}_{(\pi+\overline{\pi})^2}(4p)$$

and from the previous Lemmas, we see that

**Corollary 4.5.** *With the notations as above,*

$$K \cong \mathbb{Q}(\pi_p) \implies (-1)^g \left((4p)^g + c_1(4p)^{g-1} + \ldots + c_g\right) \cdot d(K/\mathbb{Q}) \in \mathbb{Z}^2$$

We emphasize that the factor

$$\gamma_p := (-1)^g \left((4p)^g + c_1(4p)^{g-1} + \ldots + c_g\right) \tag{12}$$

has a *uniform bound* via the Triangle Inequality that is a polynomial in $\sqrt{p}$. We will call this polynomial $\psi_g(\sqrt{p})$. One may compute that, for example, for $g = 2$

$$\gamma_p = a_2^2 - 4pa_1^2 + 4pa_2 + 4p^2$$

so that $\gamma_p \leq 128p^2$; and for $g = 3$,

$$\gamma_p = -\left((4p)^3 + (2a_2 - 6p - a_1^2)(4p)^2 + (a_2^2 - 6a_2p + 9p^2 + 2a_1a_3 - 4pa_1^2)(4p) + a_3^2 - 4pa_1a_3 + 4p^2a_1^2\right)$$

so that $\gamma_p \leq 5072p^3$. We note that these polynomials for the $\gamma_p$ are indeed quadratic in all of the $a_i$.

We now proceed to the proof of Theorem 1.1.

**Proof of Theorem 1.1.**
Let $K$ be CM field of degree $2g$ and discriminant $d := d(K/\mathbb{Q})$. We sieve the sequence

$$\mathcal{A} := (\gamma_p \cdot d)_{p \leq X}$$

with the sieving set

$$\mathcal{P} := \left\{ p \,\middle|\, z < p \leq 2z \right\}$$

with $z$ to be chosen optimally later. From Corollary 4.5, it is clear that $\Pi(A, K)(X) \leq S(\mathcal{A})$. We recall that the Square Sieve states

$$S(\mathcal{A}) \leq \frac{\#\mathcal{A}}{\#\mathcal{P}} + \max_{\substack{l,q \in \mathcal{P} \\ l \neq q}} \left| \sum_{\alpha \in \mathcal{A}} \left(\frac{\alpha}{lq}\right) \right| + \frac{2}{\#\mathcal{P}} \sum_{\alpha \in \mathcal{A}} \sum_{\substack{l \in \mathcal{P} \\ (\alpha,l) \neq 1}} 1 + \frac{1}{(\#\mathcal{P})^2} \sum_{\alpha \in \mathcal{A}} \left( \sum_{\substack{l \in \mathcal{P} \\ (\alpha,l) \neq 1}} 1 \right)^2.$$

We also recall that integration by parts yields the bounds $\sum_{p \leq X} \log p \sim X$ and $\sum_{p \leq X} (\log p)^2 \sim X \log X$, and we note that $d$, being bounded by the discriminant of $\text{char}_p(X)$, is bounded by a polynomial in $X$ that depends only on $g$. Thus,

$$\#\mathcal{A} \ll \frac{X}{\log X}; \qquad \#\mathcal{P} \asymp \frac{z}{\log z}; \qquad \sum_{\substack{l \in \mathcal{P} \\ (\alpha,l) \neq 1}} 1 \ll \log \alpha;$$

$$\sum_{\alpha \in \mathcal{A}} \log \alpha \ll \pi(X) \log d + \sum_{p \leq X} \log(\psi_g(\sqrt{p})) \ll \pi(X) \log X + \pi(X) \log X \asymp X;$$

$$\sum_{\alpha \in \mathcal{A}} (\log \alpha)^2 = \sum_{\alpha \in \mathcal{A}} \left(\log d + \log \psi_g(\sqrt{p})\right)^2$$
$$\ll_g \pi(X) \log(d)^2 + \pi(X) \log X \log d + \pi(X) \log(X)^2$$
$$\ll_g X \log X.$$

It remains to bound the character sum

$$\left| \sum_{\alpha \in \mathcal{A}} \left(\frac{\alpha}{lq}\right) \right|$$

for distinct primes $l, q \in \mathcal{P}$. We have

$$\sum_{\alpha \in \mathcal{A}} \left(\frac{\alpha}{lq}\right) = \left(\frac{d}{lq}\right) \sum_{\substack{p \leq X \\ p \nmid lqN}} \left(\frac{\gamma_p}{lq}\right) + O(\log N) + 2$$

$$= \pm \sum_{\substack{c \bmod lq \\ (c,lq)=1}} \sum_{\substack{a_1,\ldots,a_g \\ \bmod lq}} \left(\frac{\gamma_p}{lq}\right) \pi_A(X, lq; a_1, \ldots, a_g, c) + O(\log N) \qquad (13)$$

where

$$\pi_A(X, lq; a_1, \ldots, a_g, c) := \#\left\{p \leq X, p \nmid lqN \mid \mathrm{char}_p(x) \equiv x^{2g} + a_1 x^{2g-1} + \ldots + a_g x^g + ca_1 x^{g-1} \ldots + c^g \bmod lq\right\} \quad (14)$$

(We ignore the possibility that $(d, lq) \neq 1$ because we wish to bound the *maximum* value of the character sum.) Now, $\widehat{\rho}$ has open image in $\mathrm{GSp}_{2g} \widehat{\mathbb{Z}}$, by assumption; so for $z \gg_A 0$,

$$\mathrm{Gal}\left(\mathbb{Q}(A[lq])/\mathbb{Q}\right) \cong \mathrm{GSp}_{2g} \mathbb{Z}/lq\mathbb{Z}.$$

and, under the above isomorphism, specifying $\mathrm{char}_p \bmod lq$ is the same as requiring the Artin symbol $\left(\frac{\mathbb{Q}(A[lq])/\mathbb{Q}}{p}\right)$ to be contained in a certain union of conjugacy classes of $\mathrm{Gal}\left(\mathbb{Q}(A[lq])/\mathbb{Q}\right)$. Then, by the Chebotarev density theorem, for $X \gg 0$,

$$\pi_A(X, lq; a_1, \ldots, a_g, c) = \frac{\#C(lq; a_1, \ldots, a_g, c)}{\#\mathrm{GSp}_{2g} \mathbb{Z}/lq\mathbb{Z}} \pi(X) + R(X; lq; a_1, \ldots, a_g, c),$$

where

$$C(lq; a_1, \ldots, a_g, c) := \left\{h \in \mathrm{GSp}_{2g} \mathbb{Z}/lq\mathbb{Z} \mid \mathrm{char}_h(x) = x^{2g} + a_1 x^{2g-1} + \ldots + a_g x^g + ca_1 x^{g-1} \ldots + c^g\right\} \quad (15)$$

is the aforementioned union of conjugacy classes, and $R(X; lq; a_1, \ldots, a_g, c)$ is the error term, bounded variously as in Theorem 3.2. We let

$$R_{lq} := \max \left|R(X; lq; a_1, \ldots, a_g, c)\right|,$$

for notational convenience, where the maximum runs over $a_i, c \in \mathbb{Z}/lq\mathbb{Z}$. The bound (6) and the Chinese Remainder Theorem yield

$$\frac{\#C(lq; a_1, \ldots, a_g, c)}{\#\mathrm{GSp}_{2g} \mathbb{Z}/lq\mathbb{Z}} = \left(\frac{l^{2g^2}}{(l-1)(l+1)^{2g^2+g}} + Q_C(l)\right) \left(\frac{q^{2g^2}}{(q-1)(q+1)^{2g^2+g}} + Q_C(q)\right)$$
$$= f(l)f(q) + f(l)Q_C(q) + f(q)Q_C(l) + Q_c(l)Q_C(q)$$

where

$$f(l) := \frac{l^{2g^2}}{(l-1)(l+1)^{2g^2+g}}.$$

Recall that $0 \leq Q_C(l) \ll l^{-3g-1}$ (and similarly for $Q_C(q)$). Now, repeatedly using the Triangle Inequality,

$$\left|\sum_{\alpha \in \mathcal{A}} \left(\frac{\alpha}{lq}\right)\right| \leq \left|\sum_{\substack{c \bmod lq \\ (c, lq) = 1}} \sum_{a_1, \ldots, a_g \bmod lq} \left(\frac{\gamma_p}{lq}\right) \pi_A(X, lq; a_1, \ldots, a_g, c)\right| + O(\log N) + 2$$

$$\ll_N \left|\sum_c \sum_{a_1, \ldots, a_g} \left(\frac{\gamma_p}{lq}\right) \frac{\#C(lq; a_1, \ldots, a_g, c)}{\#\mathrm{GSp}_{2g} \mathbb{Z}/lq\mathbb{Z}}\right| \pi(X) + \left|\sum_c \sum_{a_1, \ldots, a_g} \left(\frac{\gamma_p}{lq}\right) R(X; lq; a_1, \ldots, a_g, c)\right|$$

$$\leq f(l)f(q) \left|\sum_c \sum_{a_1, \ldots, a_g} \left(\frac{\gamma_p}{lq}\right)\right| \pi(X) + f(l) \left|\sum_c \sum_{a_1, \ldots, a_g} \left(\frac{\gamma_p}{lq}\right) Q_C(q)\right| \pi(X)$$

$$+ f(q) \left|\sum_c \sum_{a_1, \ldots, a_g} \left(\frac{\gamma_p}{lq}\right) Q_C(l)\right| \pi(X) + \left|\sum_c \sum_{a_1, \ldots, a_g} \left(\frac{\gamma_p}{lq}\right) Q_C(l) Q_C(q)\right| \pi(X)$$

$$+ \left|\sum_c \sum_{a_1, \ldots, a_g} \left(\frac{\gamma_p}{lq}\right) R(X; lq; a_1, \ldots, a_g, c)\right|$$

so that

$$\left|\sum_{\alpha \in \mathcal{A}} \left(\frac{\alpha}{lq}\right)\right| \ll_N (lq)^{-g-1} \left|\sum_c \sum_{a_1,\ldots,a_g} \left(\frac{\gamma_p}{lq}\right)\right| \pi(X) + l^{-g-1}(lq)^{g+1}q^{-3g-1}\pi(X)$$
$$+ q^{-g-1}(lq)^{g+1}l^{-3g-1}\pi(X) + (lq)^{g+1}l^{-3g-1}q^{-3g-1}\pi(X) + R_{lq}(lq)^{g+1}$$

$$\asymp z^{-2g-2} \left|\sum_{\substack{c \bmod lq \\ (c,lq)=1}} \sum_{\substack{a_1,\ldots,a_g \\ \bmod lq}} \left(\frac{\gamma_p}{lq}\right)\right| \pi(X) + z^{-2g}\pi(X) + z^{2g+2}R_{lq} \qquad (16)$$

It remains to bound the character sum in (16). Choose $i \in \{1, \ldots, g\}$ such that $\gamma_p$ is quadratic in $a_i$. Then, by Lemma 4.4, $\gamma_p = \gamma_{i,p}^{(2)}(a_i)^2 + \gamma_{i,p}^{(1)}a_i + \gamma_{i,p}^{(0)}$, and the coefficients $\gamma_{i,p}^{(k)}$ are polynomials in the other $a_j$ and in $p$. We now break up the character sum using $a_i$,

$$\sum_{a_i} \left(\frac{\gamma_p}{lq}\right) = \#\left\{a_i \bmod lq \,\Big|\, \left(\frac{\gamma_p}{l}\right) = \left(\frac{\gamma_p}{q}\right) = 1\right\} - \#\left\{a_i \bmod lq \,\Big|\, \left(\frac{\gamma_p}{l}\right) = 1, \left(\frac{\gamma_p}{q}\right) = -1\right\}$$
$$- \#\left\{a_i \bmod lq \,\Big|\, \left(\frac{\gamma_p}{l}\right) = -1, \left(\frac{\gamma_p}{q}\right) = 1\right\} + \#\left\{a_i \bmod lq \,\Big|\, \left(\frac{\gamma_p}{l}\right) = \left(\frac{\gamma_p}{q}\right) = -1\right\}. \qquad (17)$$

These numbers are related to the number of points on certain genus-0 curves over $\mathbb{Z}/l\mathbb{Z}$ and $\mathbb{Z}/q\mathbb{Z}$, as follows. Define the projective curve $\mathcal{C}/(\mathbb{Z}/lq\mathbb{Z})$ via the affine model $\mathcal{C}^\circ$ with equation

$$y^2 = \gamma_{i,p}^{(2)}x^2 + \gamma_{i,p}^{(1)}x + \gamma_{i,p}^{(0)}.$$

and let $\mathcal{C}_l^\circ, \mathcal{C}_l$ be the reductions of $\mathcal{C}^\circ, \mathcal{C}$ modulo $l$, and similarly for $q$. Then, the number of rational points

$$\#\mathcal{C}_l^\circ(\mathbb{Z}/l\mathbb{Z}) = 2 \cdot \#\left\{a_i \bmod l \,\Big|\, \left(\frac{\gamma_p}{l}\right) = 1\right\} + \epsilon_l$$

where $\epsilon_l$ is the number of rational points $(a_i, y) \in \mathcal{C}_l^\circ(\mathbb{Z}/l\mathbb{Z})$ such that $\gamma_p \equiv 0 \bmod l$. Similarly for $q$. Now, pick a number $\xi \in \mathbb{Z}/lq\mathbb{Z}$ which is neither a square mod $l$ nor mod $q$. Then, by a similar argument, if we define the projective curve $\mathcal{C}'$ by the affine model $\mathcal{C}'^\circ$ with equation

$$y^2 = \xi\left(\gamma_{i,p}^{(2)}x^2 + \gamma_{i,p}^{(1)}x + \gamma_{i,p}^{(0)}\right)$$

and the reductions $\mathcal{C}_l'^\circ, \mathcal{C}_l'$ modulo $l$ (and similarly for $q$), then the number of rational points

$$\#\mathcal{C}'^\circ(\mathbb{Z}/l\mathbb{Z}) = 2 \cdot \#\left\{a_i \bmod l \,\Big|\, \left(\frac{\gamma_p}{l}\right) = -1\right\} + \epsilon_l'$$

with $\epsilon_l'$ defined analogously. Also denote by $\epsilon_q$, and $\epsilon_q'$ the analogous quantities for $q$. Then, by the Chinese Remainder Theorem,

$$\#\left\{a_i \bmod lq \,\Big|\, \left(\frac{\gamma_p}{l}\right) = 1 \neq \left(\frac{\gamma_p}{q}\right) = 1\right\} = \frac{1}{2}\left(\#\mathcal{C}_l^\circ(\mathbb{Z}/l\mathbb{Z}) - \epsilon_l\right) \cdot \frac{1}{2}\left(\#\mathcal{C}_q^\circ(\mathbb{Z}/q\mathbb{Z}) - \epsilon_q\right);$$

$$\#\left\{a_i \bmod lq \,\Big|\, \left(\frac{\gamma_p}{l}\right) = -1, \left(\frac{\gamma_p}{q}\right) = 1\right\} = \frac{1}{2}\left(\#\mathcal{C}_l'^\circ(\mathbb{Z}/l\mathbb{Z}) - \epsilon_l'\right) \cdot \frac{1}{2}\left(\#\mathcal{C}_q^\circ(\mathbb{Z}/q\mathbb{Z}) - \epsilon_q\right);$$

and so on for the other two terms in (17).

Assume for the moment that $\mathcal{C}_l$ is irreducible. Then, $\mathcal{C}_l$ is an irreducible genus-0 curve with a rational point. Thus, $\mathcal{C}_l \cong \mathbb{P}^1_{\mathbb{Z}/l\mathbb{Z}}$, so that $\#\mathcal{C}_l(\mathbb{Z}/l\mathbb{Z}) = l + 1$. Similarly if $\mathcal{C}_l', \mathcal{C}_q$, and $\mathcal{C}_q'$ are irreducible.

Now, $\mathcal{C}_l$ and $\mathcal{C}'_l$ are reducible iff the discriminant

$$\left(\gamma_{i,p}^{(1)}\right)^2 - 4\gamma_{i,p}^{(2)}\gamma_{i,p}^{(0)} \equiv 0 \bmod l \tag{18}$$

and similarly with $q$. Equation 18 defines a hypersurface $\mathcal{Z}_l \hookrightarrow \mathbb{A}_{\mathbb{Z}/l\mathbb{Z}}^{g-1}$ of degree at most 4, which thus has $O_g(1)$ many irreducible components. Thus, by Theorem 3.4 the number of rational points $\mathcal{Z}_l\left(\mathbb{Z}/l\mathbb{Z}\right) \ll_g l^{g-2}$. Similarly, we get a hypersurface $\mathcal{Z}_q \hookrightarrow \mathbb{A}_{\mathbb{Z}/q\mathbb{Z}}^{g-1}$ with $\ll_g q^{g-2}$ many rational points. Thus, by the Chinese Remainder Theorem, all of the curves $\mathcal{C}_l, \mathcal{C}'_l, \mathcal{C}_q$, and $\mathcal{C}'_q$ are irreducible when the numbers $(a_j)_{j\neq i} \in (\mathbb{Z}/lq\mathbb{Z})^{g-1}$ are outside a set $\mathcal{Z}$ of size $O(z^{2g-3})$.

For notational convenience, denote $\widehat{a} = (a_j)_{j\neq i} \in \left(\mathbb{Z}/lq\mathbb{Z}\right)^{g-1}$. Then, continuing from (17),

$$\left|\sum_{\widehat{a}}\sum_{a_i}\left(\frac{\gamma_p}{lq}\right)\right| \leq \left|\sum_{\widehat{a}\in\mathcal{Z}}\sum_{a_i}\left(\frac{\gamma_p}{lq}\right)\right| + \left|\sum_{\widehat{a}\notin\mathcal{Z}}\sum_{a_i}\left(\frac{\gamma_p}{lq}\right)\right|$$

$$= O_g\left(z^{2g-3}\right)\cdot(lq) + \left|\sum_{\widehat{a}\notin\mathcal{Z}}\sum_{a_i}\left(\frac{\gamma_p}{lq}\right)\right|$$

$$\ll_g z^{2g-1} + \left|\sum_{\widehat{a}\notin\mathcal{Z}}\sum_{a_i}\left(\frac{\gamma_p}{lq}\right)\right|$$

We briefly let $\delta$ (with appropriate subscripts and superscripts) denote the number of rational points at infinity of the projective curve corresponding to the subscripts and superscripts. We thus have, for $\widehat{a} \notin \mathcal{Z}$, continuing from (17),

$$\sum_{a_i}\left(\frac{\gamma_p}{lq}\right) = \frac{1}{4}\left(\#\mathcal{C}_l^\circ(\mathbb{Z}/l\mathbb{Z}) - \epsilon_l\right)\left(\#\mathcal{C}_q^\circ(\mathbb{Z}/q\mathbb{Z}) - \epsilon_q\right) + \frac{1}{4}\left(\#\mathcal{C}_l'^\circ(\mathbb{Z}/l\mathbb{Z}) - \epsilon_l'\right)\left(\#\mathcal{C}_q'^\circ(\mathbb{Z}/q\mathbb{Z}) - \epsilon_q'\right)$$

$$- \frac{1}{4}\left(\#\mathcal{C}_l'^\circ(\mathbb{Z}/l\mathbb{Z}) - \epsilon_l'\right)\left(\#\mathcal{C}_q^\circ(\mathbb{Z}/q\mathbb{Z}) - \epsilon_q\right) - \frac{1}{4}\left(\#\mathcal{C}_l^\circ(\mathbb{Z}/l\mathbb{Z}) - \epsilon_l\right)\left(\#\mathcal{C}_q'^\circ(\mathbb{Z}/q\mathbb{Z}) - \epsilon_q'\right)$$

$$= \frac{1}{4}\left(l + 1 - \delta_l - \epsilon_l\right)\left(q + 1 - \delta_q - \epsilon_q\right) + \frac{1}{4}\left(l + 1 - \delta_l' - \epsilon_l'\right)\left(q + 1 - \delta_q' - \epsilon_q'\right)$$

$$- \frac{1}{4}\left(l + 1 - \delta_l' - \epsilon_l'\right)\left(q + 1 - \delta_q - \epsilon_q\right) - \frac{1}{4}\left(l + 1 - \delta_l - \epsilon_l\right)\left(q + 1 - \delta_q' - \epsilon_q'\right)$$

$$= (\delta_l + \epsilon_l)(\delta_q + \epsilon_q) + (\delta_l' + \epsilon_l')(\delta_q' + \epsilon_q') - (\delta_l + \epsilon_l)(\delta_q' + \epsilon_q') - (\delta_l' + \epsilon_l')(\delta_q + \epsilon_q).$$

$$= O(1).$$

Thus,

$$\left|\sum_{\widehat{a}}\sum_{a_i}\left(\frac{\gamma_p}{lq}\right)\right| \ll_g z^{2g-1} + \#\left((\mathbb{Z}/lq\mathbb{Z})^{g-1} - \mathcal{Z}\right)\cdot O(1) \ll z^{2g-1},$$

so that, continuing from (16),

$$\left|\sum_{\alpha\in\mathcal{A}}\left(\frac{\alpha}{lq}\right)\right| \ll_{N,g} z^{-2g-2}z^{2g-1}\pi(X) + z^{-2g}\pi(X) + z^{2g+2}R_{lq}$$

$$\ll z^{-3}\pi(X) + z^{2g+2}R_{lq}$$

Thus, putting it together,

$$S(\mathcal{A}) \ll_N \frac{X\log z}{z\log X} + z^{-3}\pi(X) + z^{2g-1}\max_{\substack{l,q\in\mathcal{P}\\l\neq q}} R_{lq} + \frac{2\log z}{z}X + \frac{(\log z)^2}{z^2}X\log X$$

so

$$\boxed{S(\mathcal{A}) \ll \frac{X \log z}{z} + \frac{(\log z)^2 X \log X}{z^2} + z^{2g+2} \max_{l,q \in \mathcal{P}} R_{lq}(X).}$$

## 4.1 Under GRH.

Let $L_{lq} := \mathbb{Q}(A[lq])$, $n(lq) := [L_{lq} : \mathbb{Q}]$, and $d(lq) := d(L_{lq}/\mathbb{Q})$. We have the bound

$$\#C(lq; a_1, \ldots, a_g, c) \ll \left(\# \operatorname{GSp}_{2g} \mathbb{Z}/lq\mathbb{Z}\right) \cdot z^{-2g-2} \ll z^{4g^2}$$

Then, under GRH, for $X \gg 0$, Theorem 3.2 yields

$$\max R_{lq}(X) = O\left(\max_{L=L_{lq}} (\#C) X^{1/2} \left(\frac{\log|d_L|}{n_L} + \log X\right)\right)$$

$$= O\left(z^{4g^2} X^{1/2} \max_{L=L_{lq}} \left(\frac{\log|d_L|}{n_L} + \log X\right)\right)$$

where $\operatorname{Gal} L_{lq}/\mathbb{Q} = \operatorname{GSp}_{2g} \mathbb{Z}/lq\mathbb{Z}$, so $n(lq) \asymp z^{4g^2+2g+2}$. We also have, by Lemma 3.3,

$$\log|d(lq)| \leq n(lq) \log\left(\prod_{p \in \mathcal{P}(L_{lq}/\mathbb{Q})} p\right) + n(lq) \log n(lq).$$

But the only primes that ramify in $L_{lq}$ divide $lqN$. Thus,

$$\max R_{lq}(X) = O\left(z^{4g^2} X^{1/2} \left(\max \log(lqN) + \max \log(n(lq)) + \log X\right)\right)$$

$$= O_g\left(z^{4g^2} X^{1/2} \left(\log(z^2 N) + \log(z) + \log X\right)\right)$$

$$= O_{N,g}\left(z^{4g^2} X^{1/2} (\log z + \log X)\right).$$

Thus,

$$S(\mathcal{A}) \ll_{N,g} \frac{X \log z}{z} + \frac{(\log z)^2 X \log X}{z^2} + z^{2g+2} \left(z^{4g^2} X^{1/2} (\log z + \log X)\right).$$

We will choose $z$ so that $\log z \asymp \log X$. Then,

$$S(\mathcal{A}) \ll_{N,g} \frac{X \log X}{z} + z^{4g^2+2g+2} X^{1/2} \log X$$

We choose $z := X^{1/(8g^2+4g+6)}$, which yields $S(\mathcal{A}) \ll_{N,g} X^{1-1/(8g^2+4g+6)} \log X$. □

## 4.2 Under GRH + AHC.

Let the notation be as above. Under GRH and AHC, for $X \gg 0$, Theorem 3.2 yields

$$\max R_{lq}(X) = O\left(\max_{L=L_{lq}} (\#C)^{1/2} X^{1/2} \left(\log M(L/\mathbb{Q}) + \log X\right)\right)$$

$$= O_{N,g}\left(z^{2g^2} X^{1/2} (\log z + \log X)\right)$$

and thus

$$S(\mathcal{A}) \ll_{N,g} \frac{X \log z}{z} + \frac{(\log z)^2 X \log X}{z^2} + z^{2g+2} \left(z^{2g^2} X^{1/2} (\log z + \log X)\right).$$

We choose $z := X^{1/(4g^2+4g+6)}$, which yields $S(\mathcal{A}) \ll_{N,g} X^{1-1/(4g^2+4g+6)} \log X$. □

## 4.3 Under GRH + AHC + PCC.

Let the notation be as above. Under GRH, AHC, and PCC, for $X \gg 0$, Theorem 3.2 yields

$$\max R_{lq}(X) = O\left(\max(\#C)^{1/2} X^{1/2} \left(\frac{\#\widetilde{G}}{\#G}\right)^{1/4} (\log M(L/\mathbb{Q}) + \log X)\right)$$

$$= O_{N,g}\left(z^{2g^2} X^{1/2} \left(\frac{\max \#\widetilde{\mathrm{GSp}_4 \mathbb{Z}/lq\mathbb{Z}}}{z^{4g^2+2g+2}}\right)^{1/4} (\log z + \log X)\right)$$

Thus, with Lemma (3.7),

$$\max R_{lq}(X) = O_{N,g}\left(z^{2g^2} X^{1/2} \left(\frac{z^{2g+2}}{z^{4g^2+2g+2}}\right)^{1/4} (\log z + \log X)\right)$$

$$= O_{N,g}\left(z^{g^2} X^{1/2} (\log z + \log X)\right)$$

and thus

$$S(\mathcal{A}) \ll_{N,g} \frac{X \log z}{z} + \frac{(\log z)^2 X \log X}{z^2} + z^{2g+2}\left(z^{g^2} X^{1/2} (\log z + \log X)\right).$$

We choose $z := X^{1/(2g^2+4g+6)}$, which yields $S(\mathcal{A}) \ll_{N,g} X^{1-1/(2g^2+4g+6)} \log X$. □

## 4.4 Unconditionally.

Let the notation be as above. We recall part 4 of Theorem 3.2. Unconditionally, for a number field $L$, there exist constants $A, B, B' > 0$ such that when

$$\log X \geq B'(\#G)\left(\log|d_L|\right)^2,$$

we have

$$R(X) \ll \frac{\#C}{\#G} \mathrm{li}\left(X \exp\left(-B \frac{\log X}{\max\{|d_L|^{1/n_L}, \log|d_L|\}}\right)\right) + (\#\widetilde{C}) X \exp\left(-A\sqrt{\frac{\log X}{n_L}}\right). \tag{19}$$

We recall Lemma 3.3, which states

$$\frac{n_L}{2} \sum_{p \in \mathcal{P}(L/\mathbb{Q})} \log p \leq \log|d_L| \leq (n_L - 1) \sum_{p \in \mathcal{P}(L/\mathbb{Q})} \log p + n_L \log n_L.$$

Thus, with $L = \mathbb{Q}(A[lq])/\mathbb{Q}$,

$$\log|d_L| \leq z^{4g^2+2g+2}\left(\log(lqN) + \log(z^{2g^2+g+1})\right) \ll_{N,g} z^{4g^2+2g+2} \log z.$$

Now, $l$ and $q$ do ramify in $\mathbb{Q}(A[lq])/\mathbb{Q}$, since the existence of the Weil pairing on $A[lq]$ implies that $\mathbb{Q}(A[lq])/\mathbb{Q}$ contains an $(lq)^{\mathrm{th}}$ root of unity. Thus,

$$\log|d_L| \gg_N z^{4g^2+2g+2} \log(lq) \asymp z^{4g^2+2g+2} \log z$$

Also,

$$|d_L|^{1/n_L} \geq \left(\prod_{p \in \mathcal{P}(L/\mathbb{Q})} p\right)^{1/2} \geq (lq)^{1/2} \asymp z$$

and
$$|d_L|^{1/n_L} \leq n_L \prod_{p \in \mathcal{P}(L/\mathbb{Q})} p \leq z^{4g^2+2g+2}(lqN) \ll_N z^{4g^2+2g+4}.$$

Thus, the requirement
$$\log X \geq B'(\#G)\left(\log|d_L|\right)^2 \asymp_{N,g} B' z^{8g^2+6g+4}(\log z)^2 \tag{20}$$

is satisfied with the choice
$$z := c' \frac{(\log X)^{1/(8g^2+6g+4)}}{(\log\log X)^{1/(4g^2+3g+2)}}$$

for a certain positive constant $c'$ depending only on $N$ and $g$. The reader may check that there exists such a $c'$ so that (20) is satisfied with this choice of $z$. Moreover, we see from the above that $\max\{|d_L|^{1/n_L}, \log|d_L|\} \ll_N z^{4g^2+2g+4}$.

For $l, q \in \mathcal{P}$, arguments above show that
$$\frac{\#C}{\#\operatorname{GSp}_4\mathbb{Z}/lq\mathbb{Z}} \ll z^{-2g-2}.$$

Using the approximation $\operatorname{li} t \sim \frac{t}{\log t}$, we then have

$$\frac{\#C}{\#G} \operatorname{li}\left(X \exp\left(-B\frac{\log X}{\max\{|d_L|^{1/n_L}, \log|d_L|\}}\right)\right) \ll z^{-2g-2} \frac{\left(X \exp\left(-B\frac{\log X}{\max\{|d_L|^{1/n_L}, \log|d_L|\}}\right)\right)}{\log\left(X \exp\left(-B\frac{\log X}{\max\{|d_L|^{1/n_L}, \log|d_L|\}}\right)\right)}$$

$$\ll_N \frac{z^{-2g-2} X \exp\left(-B\frac{\log X}{z^{4g^2+2g+4}\log z}\right)}{\log X - B(\log X)z^{-(4g^2+2g+4)}(\log z)^{-1}}$$

$$\ll \frac{X^{1-Bz^{-(4g^2+2g+4)}(\log z)^{-1}}}{z^{4g^2+2g+4}\log X}$$

From our choice of $z$, (24), the bounds above, and the weak bound $\#\widetilde{C} \leq \#\widetilde{\operatorname{GSp}_{2g}\mathbb{Z}/lq\mathbb{Z}} \asymp z^{2g+2}$, we obtain (after a calculation which we omit; see Section 4 of [CFM05]) the bounds

$$\max_{\substack{l,q\in\mathcal{P}\\l\neq q}}\left|\sum_{\alpha\in\mathcal{A}}\left(\frac{\alpha}{lq}\right)\right| \ll_N \frac{X}{z\log X};$$

$$\sum_{\alpha\in\mathcal{A}}\sum_{\substack{l\in\mathcal{P}\\(\alpha,l)\neq 1}} 1 \ll_N \frac{X}{\log X}\nu_z(d);$$

$$\sum_{\alpha\in\mathcal{A}}\left(\sum_{\substack{l\in\mathcal{P}\\(\alpha,l)\neq 1}} 1\right)^2 \ll_N \frac{X}{\log X}\left(\nu_z(d) + (\nu_z(d))^2\right);$$

where $\nu_z(d)$ is the number of distinct prime divisors of $d$ less than or equal to $z$. Thus, from the Square Sieve and the trivial bounds $\nu_z(d) \leq \nu(d)$ and $\nu_z(d) \leq \pi(z)$ we obtain

$$S(\mathcal{A}) \ll_{N,g} \frac{X\log z}{z\log X}(1+\nu_z(d)) \ll_{N,g} \frac{X(\log\log X)^{1+1/(4g^2+3g+2)}}{(\log X)^{1+1/(8g^2+6g+4)}}(1+\nu_z(\mathrm{d}(K/\mathbb{Q}))). \tag{21}$$

$\square$

# 5 Proof of Theorem 1.2

In this Section, $A/\mathbb{Q}$ is a principally polarized abelian surface with $\mathrm{End}\left(A_{\overline{\mathbb{Q}}}\right) \cong \mathbb{Z}$. This implies that $A$ is simple. As mentioned in Remark 1.3, works of Serre show that its adelic Galois representation $\widehat{\rho}$ has open image in $\mathrm{GSp}_4\widehat{\mathbb{Z}}$. Let $N$ be the conductor of $A$, and let $F/\mathbb{Q}$ be a real quadratic number field.

Let $p \nmid N$ be a prime of good, ordinary, non-split reduction for $A$. Then, by Honda-Tate theory, the endomorphism algebra $K := \mathrm{End}(A_p) \otimes \mathbb{Q} = \mathbb{Q}(\pi_p)$ is a quartic CM field. Let $K_0$ be the totally real quadratic subfield of $K$. Then,

$$K_0 = \mathbb{Q}(\sqrt{d}), \quad K = \mathbb{Q}(\pi_p) = K_0(\sqrt{r})$$

for some squarefree rational integer $d > 0$, and some totally negative integer $r \in \mathcal{O}_{K_0}$.

As in Section 4, the characteristic polynomial of the Frobenius endomorphism $\pi_p$ has the shape

$$\mathrm{char}_p(x) = x^4 + a_{1,p}x^3 + a_{2,p}x^2 + pa_{1,p}x + p^2,$$

and the Triangle Inequality yields

$$|a_{1,p}| \leq \binom{4}{1}\sqrt{p} = 4\sqrt{p}; \qquad |a_{2,p}| \leq \binom{4}{2}(\sqrt{p})^2 = 6p. \tag{22}$$

**Remark 5.1.** *Since $A$ is a simple abelian surface, $A$ is the Jacobian of some smooth curve $C$ of genus 2; it is well known that $a_{1,p}$ and $a_{2,p}$ may be expressed in terms of the number of $\mathbb{F}_p$- and $\mathbb{F}_{p^2}$-points of the reduction of $C$ mod $p$, as one has the formula of Hasse-Weil,*

$$\#C_p(\mathbb{F}_{p^k}) = p^k + 1 - \sum \lambda^k$$

*where the sum is over the roots $\lambda \in \overline{\mathbb{Q}}$ of $\mathrm{char}_p$. Letting $N_k := \#C_p(\mathbb{F}_{p^k})$, this yields the formulas*

$$a_1 = p + 1 - N_1; \qquad a_2 = \frac{1}{2}\left(N_2 + N_1(N_1 - 2p - 2)\right).$$

Now, the following lemma allows us to apply the Square Sieve to $\Pi(A, F)$.

**Lemma 5.2.** *Let the notation be as above, but with $d > 0$ an arbitrary squarefree rational integer. Then,*

$$K_0 \cong \mathbb{Q}(\sqrt{d}) \iff d(a_1^2 - 4a_2 + 8p) \text{ is a square.}$$

*Proof.* Let $x \mapsto \overline{x}$ be the complex conjugation of $K/K_0$. Then,

$$K_0 = \mathbb{Q}(\pi + \overline{\pi}) = \mathbb{Q}\left(\pi + \frac{p}{\pi}\right).$$

(Note that $\pi + \overline{\pi} \notin \mathbb{Q}$ because $\pi$ satisfies $x^2 - (\pi + \overline{\pi})x + p = 0$ and $[\mathbb{Q}(\pi) : \mathbb{Q}] = 4$.) Let $\beta = \pi + p/\pi$. Then, for $m, n \in \mathbb{Z}$,

$$\beta^2 + m\beta + n = 0 \iff \pi^4 + m\pi^3 + (2p + n)\pi^2 + pm\pi + p^2 = 0$$

so that the minimal polynomial of $\beta$ is $x^2 + a_1x + a_2 - 2p$. The result follows from the requirement that $d$ and the discriminant $a_1^2 - 4(a_2 - 2p)$ must have the same squarefree part. $\square$

We now proceed to the proof of Theorem 1.2. Because of its similarity to the proof of Theorem 1.1, we abbreviate some parts of the proof.

**Proof of Theorem 1.2.**
We apply the Square Sieve to the sequence

$$\mathcal{A} := \left(d(a_{1,p}^2 - 4a_{2,p} + 8p)\right)_{p \leq X}$$

with the sieving set

$$\mathcal{P} := \left\{ p \mid z < p \leq 2z \right\}$$

with $z$ to be chosen optimally later. From Lemma 5.2, it is clear that $\Pi(A, F)(X) \leq S(\mathcal{A})$.

We recall that the Square Sieve states

$$S(\mathcal{A}) \leq \frac{\#\mathcal{A}}{\#\mathcal{P}} + \max_{\substack{l,q \in \mathcal{P} \\ l \neq q}} \left| \sum_{\alpha \in \mathcal{A}} \left( \frac{\alpha}{lq} \right) \right| + \frac{2}{\#\mathcal{P}} \sum_{\alpha \in \mathcal{A}} \sum_{\substack{l \in \mathcal{P} \\ (\alpha, l) \neq 1}} 1 + \frac{1}{(\#\mathcal{P})^2} \sum_{\alpha \in \mathcal{A}} \left( \sum_{\substack{l \in \mathcal{P} \\ (\alpha, l) \neq 1}} 1 \right)^2.$$

We have again the bounds

$$\#\mathcal{A} \ll \frac{X}{\log X}; \qquad \#\mathcal{P} \asymp \frac{z}{\log z}; \qquad \sum_{\substack{l \in \mathcal{P} \\ (\alpha, l) \neq 1}} 1 \ll \log \alpha;$$

$$\sum_{\alpha \in \mathcal{A}} \log \alpha \ll X; \qquad \sum_{\alpha \in \mathcal{A}} (\log \alpha)^2 \ll X \log X.$$

It remains to bound the character sum

$$\left| \sum_{\alpha \in \mathcal{A}} \left( \frac{\alpha}{lq} \right) \right|$$

for distinct primes $l, q \in \mathcal{P}$. We have

$$\sum_{\alpha \in \mathcal{A}} \left( \frac{\alpha}{lq} \right) = \left( \frac{d}{lq} \right) \sum_{\substack{p \leq X \\ p \nmid lqN}} \left( \frac{a_{1,p}^2 - 4a_{2,p} + 8p}{lq} \right) + O(\log N)$$

$$= \pm \sum_{\substack{c \bmod lq \\ (c, lq) = 1}} \sum_{\substack{a_1, a_2 \\ \bmod lq}} \left( \frac{a_1^2 - 4a_2 + 8c}{lq} \right) \pi_A(X, lq; a_1, a_2, c) + O(\log N) \qquad (23)$$

where $\pi_A(X, lq; a_1, a_2, c)$ is defined as in (14). Then, by the Chebotarev density theorem, for $X \gg 0$,

$$\pi_A(X, lq; a_1, a_2, c) = \frac{\#C(lq; a_1, a_2, c)}{\#\operatorname{GSp}_4 \mathbb{Z}/lq\mathbb{Z}} \pi(X) + R(X; lq; a_1, a_2, c),$$

where $C(lq; a_1, a_2, c)$ is defined as in (15), and $R(X; lq; a_1, a_2, c)$ is the error term, bounded variously as in Theorem 3.2. We let

$$R_{lq} := \max \left| R(X; lq; a_1, a_2, c) \right|$$

for notational convenience, where the max runs over $a_1, a_2, x \in \mathbb{Z}/lq\mathbb{Z}$.

The bound (6) with $g = 2$ and the Chinese Remainder Theorem yields

$$\frac{\#C(lq; a_1, a_2, c)}{\#\operatorname{GSp}_4 \mathbb{Z}/lq\mathbb{Z}} = \frac{l^8}{(l-1)(l+1)^{10}} \cdot \frac{q^8}{(q-1)(q+1)^{10}} + Q(lq; a_1, a_2, c)$$

where the error term satisfies $0 \leq Q(lq; a_1, a_2, c) = l^{-3} \cdot O(q^{-7}) + q^{-3} \cdot O(l^{-3}) + O\left((lq)^{-7}\right) = O\left(z^{-10}\right)$. Thus,

$$\left| \sum_{\alpha \in \mathcal{A}} \left( \frac{\alpha}{lq} \right) \right| \ll \left| \sum_{\substack{c \bmod lq \\ (c, lq) = 1}} \sum_{\substack{a_1, a_2 \\ \bmod lq}} \left( \frac{a_1^2 - 4a_2 + 8c}{lq} \right) \left( \frac{l^8 q^8 \pi(X)}{(l-1)(l+1)^{10}(q-1)(q+1)^{10}} + Q(lq; a_1, a_2, c)\pi(X) \right) \right|$$

$$+ (lq)^3 R_{lq}(X) + O(\log N),$$

Now, by the orthogonality of characters, once $l, q > 2$,

$$\sum_{a_2 \bmod lq} \left(\frac{a_1^2 - 4a_2 + 8c}{lq}\right) = 0$$

Thus,

$$\left|\sum_{\alpha \in \mathcal{A}} \left(\frac{\alpha}{lq}\right)\right| \ll_N \left|\sum_{\substack{c \bmod lq \\ (c, lq) = 1}} \sum_{\substack{a_1, a_2 \\ \bmod lq}} \left(\frac{a_1^2 - 4a_2 + 8c}{lq}\right) Q(lq; a_1, a_2, c)\right| \pi(X) + (lq)^3 R_{lq}(X),$$

and thus, by a similar argument as what led to (16),

$$\left|\sum_{\alpha \in \mathcal{A}} \left(\frac{\alpha}{lq}\right)\right| \ll_N (lq)^3 \left(l^{-3} q^{-7} + l^{-7} q^{-3}\right) \pi(X) + (lq)^3 R_{lq}(X)$$

$$\asymp z^{-4} \frac{X}{\log X} + z^6 R_{lq}(X).$$

Putting it together,

$$S(\mathcal{A}) \ll_N \frac{X \log z}{z \log X} + \left(\frac{X}{z^4 \log X} + z^6 \max_{\substack{l, q \in \mathcal{P} \\ l \neq q}} R_{lq}(X)\right) + \frac{2 \log z}{z} X + \frac{(\log z)^2}{z^2} X \log X$$

and thus

$$\boxed{S(\mathcal{A}) \ll_N \frac{X \log z}{z} + \frac{(\log z)^2 X \log X}{z^2} + z^6 \max R_{lq}(X).}$$

## 5.1 Under GRH.

Let $L = L_{lq} := \mathbb{Q}(A[lq])$, $n(lq) := [L_{lq} : \mathbb{Q}]$, and $d(lq) := d(L_{lq}/\mathbb{Q})$. We have the bound

$$\#C(lq; a_1, a_2, c) \ll (\#\operatorname{GSp}_4 \mathbb{Z}/lq\mathbb{Z}) \frac{(lq)^8}{(l-1)(l+1)^{10}(q-1)(q+1)^{10}} \ll (lq)^8$$

Then, under GRH, for $X \gg 0$, Theorem 3.2 yields

$$\max R_{lq}(X) = O\left(\max_{L = L_{lq}} (\#C) X^{1/2} \left(\frac{\log|d_L|}{n_L} + \log X\right)\right)$$

$$= O\left(z^{16} X^{1/2} \max_{L = L_{lq}} \left(\frac{\log|d_L|}{n_L} + \log X\right)\right)$$

where $\operatorname{Gal} L_{lq}/\mathbb{Q} = \operatorname{GSp}_4 \mathbb{Z}/lq\mathbb{Z}$, so $n(lq) \asymp (lq)^{10} \asymp z^{20}$. We also have, by Lemma (3.3),

$$\log|d(lq)| \leq n(lq) \log\left(\prod_{p \in \mathcal{P}(L_{lq}/\mathbb{Q})} p\right) + n(lq) \log n(lq).$$

But the only primes that ramify in $L_{lq}$ divide $lqN$. So,

$$\max R_{lq}(X) = O\left(z^{16} X^{1/2} \left(\max \log(lqN) + \max \log(n_L) + \log X\right)\right)$$

so

$$\max R_{lq}(X) = O\left(z^{16} X^{1/2} \left(\log(z^2 N) + \log(z^{10} z^{10}) + \log X\right)\right)$$
$$= O_N\left(z^{16} X^{1/2} \left(\log z + \log X\right)\right).$$

Thus,

$$S(\mathcal{A}) \ll_N \frac{X \log z}{z} + \frac{(\log z)^2 X \log X}{z^2} + z^6 \left(z^{16} X^{1/2} \left(\log z + \log X\right)\right).$$

We choose $z := X^{1/46}$, which yields $S(\mathcal{A}) \ll_N X^{45/46} \log X$. □

## 5.2 Under GRH + AHC.

Let the notation be as above. Under GRH and AHC, for $X \gg 0$, Theorem 3.2 yields

$$\max R_{lq}(X) = O\left(\max_{L=L_{lq}} (\#C)^{1/2} X^{1/2} \left(\log M(L/\mathbb{Q}) + \log X\right)\right)$$
$$= O\left(z^8 X^{1/2} \left(\log(z^2 N) + \log X\right)\right)$$
$$= O_N\left(z^8 X^{1/2} \left(\log z + \log X\right)\right)$$

and thus

$$S(\mathcal{A}) \ll_N \frac{X \log z}{z} + \frac{(\log z)^2 X \log X}{z^2} + z^6 \left(z^8 X^{1/2} \left(\log z + \log X\right)\right).$$

We choose $z := X^{1/30}$, which yields $S(\mathcal{A}) \ll_N X^{29/30} \log X$. □

## 5.3 Under GRH + AHC + PCC.

Let the notation be as above. Under GRH, AHC, and PCC, for $X \gg 0$, Theorem 3.2 yields

$$\max R_{lq}(X) = O\left(\max(\#C)^{1/2} X^{1/2} \left(\frac{\#\widetilde{G}}{\#G}\right)^{1/4} \left(\log M(L/\mathbb{Q}) + \log X\right)\right)$$
$$= O\left(z^8 X^{1/2} \left(\frac{\max \#\widetilde{\mathrm{GSp}_4 \mathbb{Z}/lq\mathbb{Z}}}{z^{20}}\right)^{1/4} \left(\log(z^2 N) + \log X\right)\right)$$

Thus, with Lemma (3.7),

$$\max R_{lq}(X) = O\left(z^8 X^{1/2} \left(\frac{z^6}{z^{20}}\right)^{1/4} \left(\log(z^2 N) + \log X\right)\right)$$
$$= O_N\left(z^{9/2} X^{1/2} (\log z + \log X)\right)$$

and thus

$$S(\mathcal{A}) \ll_N \frac{X \log z}{z} + \frac{(\log z)^2 X \log X}{z^2} + z^6 \left(z^{9/2} X^{1/2} \left(\log z + \log X\right)\right).$$

We choose $z := X^{1/23}$, which yields $S(\mathcal{A}) \ll_N X^{22/23} \log X$. □

## 5.4 Unconditionally.

Let the notation be as above. We recall part 4 of Theorem 3.2. Unconditionally, there exist constants $A, B, B' > 0$ such that when

$$\log X \geq B'(\#G)\left(\log|d_L|\right)^2,$$

we have

$$R(X) \ll \frac{\#C}{\#G} \operatorname{li}\left(X \exp\left(-B\frac{\log X}{\max\{|d_L|^{1/n_L}, \log|d_L|\}}\right)\right) + (\#\tilde{C})X \exp\left(-A\sqrt{\frac{\log X}{n_L}}\right). \tag{24}$$

We recall Lemma (3.3), which states

$$\frac{n_L}{2} \sum_{p \in \mathcal{P}(L/\mathbb{Q})} \log p \leq \log|d_L| \leq (n_L - 1) \sum_{p \in \mathcal{P}(L/\mathbb{Q})} \log p + n_L \log n_L.$$

Thus, by arguments identical as in Subsection 4.4,

$$\log|d_L| \ll_N z^{22} \log z; \qquad \log|d_L| \gg_N z^{22} \log z;$$

$$|d_L|^{1/n_L} \gg z; \qquad |d_L|^{1/n_L} \ll_N z^{24}$$

so that $\max\{|d_L|^{1/n_L}, \log|d_L|\} \ll_N z^{24}$. We see that the requirement

$$\log X \geq B'(\#G)\left(\log|d_L|\right)^2 \asymp_N B'z^{66}(\log z)^2 \tag{25}$$

is satisfied with the choice

$$z := c'\frac{(\log X)^{1/66}}{(\log\log X)^{1/33}}$$

for a certain positive constant $c'$ depending only on $N$. The reader may check that there exists such a $c'$ so that 25 is satisfied with this choice of $z$.

For $l, q \in \mathcal{P}$, arguments above show that

$$\frac{\#C}{\#\operatorname{GSp}_4 \mathbb{Z}/lq\mathbb{Z}} \ll z^{-6}.$$

Using the approximation $\operatorname{li} t \sim \frac{t}{\log t}$, we then have

$$\frac{\#C}{\#G} \operatorname{li}\left(X \exp\left(-B\frac{\log X}{\max\{|d_L|^{1/n_L}, \log|d_L|\}}\right)\right) \ll z^{-6}\frac{\left(X \exp\left(-B\frac{\log X}{\max\{|d_L|^{1/n_L}, \log|d_L|\}}\right)\right)}{\log\left(X \exp\left(-B\frac{\log X}{\max\{|d_L|^{1/n_L}, \log|d_L|\}}\right)\right)}$$

$$\ll_N \frac{z^{-6} X \exp\left(-B\frac{\log X}{z^{22}\log z}\right)}{\log X - B(\log X)z^{-1/2}}$$

$$\ll \frac{X^{1 - Bz^{-22}(\log z)^{-1}}}{z^6 \log X}$$

From our choice of $z$, the bound of (24), the bounds above, and the weak bound $\#\widetilde{C} \leq \#\widetilde{\mathrm{GSp}_4\mathbb{Z}/lq\mathbb{Z}} \asymp z^6$, we obtain (after another calculation that we omit; see Section 4 of [CFM05]) the bounds

$$\max_{\substack{l,q\in\mathcal{P}\\ l\neq q}} \left|\sum_{\alpha\in\mathcal{A}} \left(\frac{\alpha}{lq}\right)\right| \ll_N \frac{X}{z^8 \log X};$$

$$\sum_{\alpha\in\mathcal{A}} \sum_{\substack{l\in\mathcal{P}\\ (\alpha,l)\neq 1}} 1 \ll_N \frac{X}{\log X}\nu_z(d);$$

$$\sum_{\alpha\in\mathcal{A}} \left(\sum_{\substack{l\in\mathcal{P}\\ (\alpha,l)\neq 1}} 1\right)^2 \ll_N \frac{X}{\log X}\left(\nu_z(d) + \nu_z(d)^2\right);$$

where $\nu_z(d)$ is the number of distinct prime divisors of $d$ less than or equal to $z$. Thus, from the Square Sieve and the trivial bounds $\nu_z(d) \leq \nu(d)$ and $\nu_z(d) \leq \pi(z)$ we obtain

$$S(\mathcal{A}) \ll_N \frac{X \log z}{z \log X}(1 + \nu(d)) \ll_N \frac{X(\log \log X)^{23/22}}{(\log X)^{67/66}}(1 + \nu(d)). \tag{26}$$

$\square$

## 6 Proof of Corollaries 1.4, 1.5, and 1.6

The proofs for Corollaries 1.4 and 1.5 are nearly identical, so for brevity we only prove the former. We mimic the argument based on the Pigeonhole Principle in [CFM05].

We recall that, if $p$ is a good ordinary non-split prime for $A$, then $d(\mathbb{Q}(\pi_p)/\mathbb{Q})$ has squarefree part dividing the number $\gamma_p$ defined in (12). Moreover, the functions $\psi_g(\sqrt{X})$ were defined precisely so that $|\gamma_p| \leq \psi_g(\sqrt{X})$ when $p \leq X$. In Section 1, we defined the field-counting function,

$$\mathcal{D}_A(X) := \left\{K \in \mathcal{D}_A \mid \mathrm{sf}(d(K/\mathbb{Q})) \leq \psi_g(\sqrt{X})\right\}$$

so that if $p$ is good ordinary non-split for $A$, then $p \leq X$ implies $\mathbb{Q}(\pi_p) \in \mathcal{D}_A(X)$.

Now, note that because $A$ has *trivial* geometric endomorphism algebra, the set of non-split primes for $A$ has density zero (see Subsection 2.4). Thus, assuming that the set of ordinary primes for $A$ has positive density $\delta$, we may write

$$\pi(X) = (1 - \delta)\pi(X) + o(\pi(X)) + \sum_{K\in\mathcal{D}_A(\infty)} \Pi(A, K)(X)$$

$$= (1 - \delta + o(1))\pi(X) + \sum_{K\in\mathcal{D}_A(X)} \Pi(A, K)(X)$$

and thus obtain

$$\#\mathcal{D}_A(X) \geq \frac{(\delta - o(1))\pi(X)}{\max_{K\in\mathcal{D}_A(X)} \Pi(A, K)} \tag{27}$$

Plugging in the various conditional asymptotic upper bounds of Theorem 1.2 on $\Pi(A, K)(X)$ yields the conditional asymptotic lower bounds of Corollary 1.4.

Unfortunately, the dependency in $d(K/\mathbb{Q})$ of the unconditional bound for $\Pi(A, K)(X)$ keeps this argument from working in the unconditional case. But to prove Corollary 1.6, we argue as follows. By Theorem 1.1, we know that each set $\Pi(A, K)$ has density zero in the set of rational primes. Yet the set of primes at which $A$ has good, ordinary, non-split reduction is assumed to have positive density. Thus, there must be infinitely many CM fields $K$ for which $\Pi(A, K) \neq \emptyset$. $\square$

# 7 Concluding Remarks and Further Directions.

We acknowledge that the technique that we use here is no longer the state of the art for the "fixed-field" question of Lang-Trotter type that pertains to elliptic curves, and we hope that the methods of other works (such as those mentioned in Subsection 2.3) may be used to improve our result.

**Remark 7.1.** *It would be interesting to consider Question 0.3 for abelian varieties other than those whose adelic Galois representation $\widehat{\rho}$ has open image in $\mathrm{GSp}_{2g}\widehat{\mathbb{Z}}$. For the methods of this paper to work, one would need to require that the image of $\widehat{\rho}$ be open in $\mathcal{G}(\widehat{\mathbb{Z}})$ for some sub-group-scheme $\mathcal{G} \hookrightarrow \mathrm{GSp}_{2g}$; in particular, one would need to require that the images $\rho_{lq}(G_{\mathbb{Q}})$ eventually "stabilize" as $\mathcal{G}(\mathbb{Z}/lq\mathbb{Z})$ as $l,q \to \infty$.*

**Remark 7.2.** *Question 0.4 and variants thereof are easily extended to non-CM abelian varieties of any dimension $g \geq 2$. Namely, we may ask about the set of primes $p$ for which any specified ring $R$ embeds into the endomorphism algebra $\mathrm{End}(A_p) \otimes \mathbb{Q}$. (Of course, one would specify that $R$ must be a ring which embeds into the endomorphism algebra of some abelian variety of dimension $g$ over a finite field.)*

One may ask questions of Lang-Trotter type about abelian varieties in analogy to other questions asked about the reductions of elliptic curves. The only such question that seems to have been studied thus far is a generalization of the "fixed-trace" question, which is studied in [Coj+16] We give a few questions below:

**Question 7.3.** *How often is the order of the group of rational points, $\#A_p(\mathbb{F}_p)$ prime, nearly-prime, or pseudoprime (to a fixed base)? On the analogous questions for elliptic curves, see for instance [Kob88; MM01; CLS09], and see [BCD11] for the study of the primality of $\#E_p(\mathbb{F}_p)$ on average.*

**Question 7.4.** *Let $F/\mathbb{Q}$ be a totally real field of degree $g$, and $K/\mathbb{Q}$ a primitive CM field of degree $2g$. What are the values of $\Pi(A,F)(X)$ and $\Pi(A,K)(X)$ on average for generic $A$? One would need to specify how to average. For $g=2$ or $g=3$, one could averaging over boxes for the coefficients of genus-$g$ curves $C/\mathbb{Q}$, considering these counting functions for the Jacobian of $C$. (Once $g \geq 4$, not all abelian varieties are Jacobians of curves.) See, for instance, [Shp13] on the analogous question for elliptic curves.*

**Question 7.5.** *Let $A_1$ and $A_2$ be abelian varieties over $\mathbb{Q}$ (generic or otherwise). How can we describe the set of primes at which both the $A_i$ are good ordinary non-split, and*

1. *$\mathbb{Q}(\pi_{p,A_1}) \cong \mathbb{Q}(\pi_{p,A_2})$? or $\mathbb{Q}(\pi_{p,A_1}) \cong \mathbb{Q}(\pi_{p,A_2}) \cong K$ for a specified primitive CM field $K$?*

2. *$\mathbb{Q}(\pi_{p,A_1})_0 \cong \mathbb{Q}(\pi_{p,A_2})_0$? or $\mathbb{Q}(\pi_{p,A_1})_0 \cong \mathbb{Q}(\pi_{p,A_2})_0 \cong F$ for a specified totally real field $F$?*

3. *$a_{i,p,A_1} = a_{i,p,A_2}$ for specified $i$? or $a_{i,p,A_1} = a_{i,p,A_2} = t$ for specified $i$ and $t$?*

4. *their characteristic polynomials of Frobenius at $p$ are equal?*

*One may also ask these questions without the requirement that $p$ be ordinary for the $A_i$.*

We hesitate to give precise conjectures on the asymptotic growth (or boundedness) of the functions $\Pi(A,K)(X)$ and $\Pi(A,F)(X)$ until after further study and collection of data. We have run small-scale experiments; as an example, let $C/\mathbb{Q}$ be the curve of genus 2 with affine model $y^2 = x^5 - 3x^4 + 2x^3 + 1$, and let $A/\mathbb{Q}$ be the Jacobian of $C$. ($C$ is the curve 3680.a.29440.1 of [LMFDB].) Since $A$ is an abelian surface without extra endomorphisms, its adelic Galois representation has open image in $\mathrm{GSp}_4\widehat{\mathbb{Z}}$, so our results apply to $A$. We found via a simple program written in Sage [SageMath] that, in fact, $\Pi(A,K)(10^6) \leq 1$ for all $K$.